\documentclass[11pt]{article}
\usepackage[margin=1in]{geometry}
\usepackage{amsmath,amssymb,amsthm}
\usepackage{xurl}
\usepackage{longtable}
\usepackage[hidelinks,breaklinks=true]{hyperref}
\numberwithin{equation}{section}
\newtheorem{theorem}{Theorem}[section]
\newtheorem{proposition}[theorem]{Proposition}
\newtheorem{lemma}[theorem]{Lemma}
\newtheorem{corollary}[theorem]{Corollary}
\newtheorem{mainthm}{Theorem}

\theoremstyle{definition}

\theoremstyle{remark}
\newtheorem{remark}[theorem]{Remark}
\newcommand{\osp}{\mathfrak{osp}}
\newcommand{\ad}{\operatorname{ad}}
\title{Inhomogeneous deformations of $\osp(N|2n)$: triviality for any number of odd
parameters}
\author{Hisashi Aoi\thanks{Department of Mathematical Sciences, Ritsumeikan University}}
\date{\today}
\begin{document}
\maketitle

\begin{abstract}
We treat a symmetrized mixed-oscillator deformation family of the
orthosymplectic Lie superalgebras $\osp(N|2n)$, $N,n\geq1$, with even mixed
coefficients $\beta_{u,i}$ and one odd square-zero parameter. We derive the deformed
bracket from an oscillator realization, show that its coefficient is determined by the
realization, and exhibit an explicit odd primitive, normalized by its vanishing on
$\mathfrak{sp}(2n)$. Every odd primitive differs from it by a unique inner derivation.
For $N=1$ this is the case $\osp(1|2n)$ treated earlier. With multiple odd parameters
the deformed bracket is the sum of the single-parameter deformations and is again
trivial. The trivialization, however, is not additive: it is the product of the
single-parameter trivializations followed by a nonsymplectic change of the bosonic
generators. The main results are formalized in Lean~4, with the algebras constructed by
generators and relations; what is not formalized is stated in the appendix.
\end{abstract}

\section{Introduction}\label{sec:introduction}

In~\cite{BS17}, Bakalov and Sullivan introduce supersymmetric bilinear forms
$(\cdot|\cdot)$ on a superspace $V=V_{\bar0}\oplus V_{\bar1}$ that they call
\emph{inhomogeneous} when $(V_{\bar0}|V_{\bar1})\neq0$ and $(V_i|V_i)\neq0$ for at
least one $i\in\{\bar0,\bar1\}$. From a skew-supersymmetric pre-oscillator form they
construct a central extension $V\oplus\mathbb CK\oplus\mathbb C\kappa$, with an even
central element $K$ and an odd central element $\kappa$. The bracket of two generators is
valued in $K$ when their parities agree, and in $\kappa$ when they differ~\cite[Section~4]{BS17}. In the relations \eqref{eq:source} below,
$[b_u,b_v]=J_{uv}$ and $\{a_i,a_j\}=g_{ij}$ are the homogeneous Weyl and Clifford
relations, and $[b_u,a_i]=\beta_{u,i}\kappa$ is the odd central-valued term added by
the mixed block $(b_u|a_i)\neq0$.

In~\cite{Aoi26} we treated $\osp(1|2n)$, with one odd generator $a$ and mixed
coefficients $\beta_u$. Here we take $N$ odd generators $a_i$, with a nondegenerate
symmetric form $g$, and a matrix of mixed coefficients $\beta_{u,i}$. Precisely, let $P$ be
the polynomial ring in the $\beta_{u,i}$ and $R=P\oplus\kappa P$ with $\kappa$ odd and
$\kappa^2=0$. Let $A_\beta$ be the associative $R$-superalgebra generated by even $b_u$ ($u\in I$) and odd
$a_i$ ($i\in K_N$) subject to \eqref{eq:source}, with $J$ the standard symplectic form, and let $A_0$ be the
algebra with $\beta=0$ (Section~\ref{sec:setting}). The symmetrized quadratic elements of the
generators, the \emph{lifts}, correspond to a basis of $\mathfrak g=\osp(N|2n)$
(Proposition~\ref{prop:osp}). Their supercommutators in $A_0$ give the bracket
$[\cdot,\cdot]_0$ of $\mathfrak g$, and those in $A_\beta$ give a bracket
$[\cdot,\cdot]_\beta=[\cdot,\cdot]_0+\kappa\Gamma_\beta$ on $\mathfrak g_R=R\otimes\mathfrak g$
(Section~\ref{sec:recovery}). The algebra $\osp(N|2n)$ is $B(m,n)$ for $N=2m+1$, $C(n+1)$ for
$N=2$, and $D(m,n)$ for $N=2m\geq4$. We prove that the resulting odd deformation is
trivial for every $N$ and $n$, by an explicit primitive.

A vanishing statement alone would not give this. The simple finite-dimensional Lie
superalgebras over $\mathbb C$ admitting deformations with odd parameters are
classified in~\cite[Theorem~7.1]{Lei23}, and the orthosymplectic ones are not among
them. The present result differs from such a statement in three ways. It is
stated over $R$, with the $\beta_{u,i}$ indeterminate. The
coefficient $\Gamma_\beta$ is recovered from the realization and shown to be
determined by it. Finally, the primitive is explicit: it is linear in $\beta$, and it is
unique once it is required to vanish on $\mathfrak{sp}(2n)$.

Our main result concerns any number of odd parameters. In Section~\ref{sec:multi} the parameter
$\kappa$ is replaced by odd $\kappa_1,\ldots,\kappa_r$, the coefficients become
$R_r=P_r\otimes\Lambda(\kappa_1,\ldots,\kappa_r)$, with $P_r$ the polynomial ring in the
$\beta^{(s)}_{u,i}$, and the mixed relation becomes
\[
 [b_u,a_i]=\gamma_{u,i},\qquad \gamma_{u,i}=\sum_{s=1}^r\beta^{(s)}_{u,i}\kappa_s .
\]

\begin{mainthm}\label{thm:A}
Let $N,n,r\geq1$. The supercommutator of two lifts in $A_\beta$ is an $R_r$-linear combination of lifts, with no
scalar term, and the resulting bracket on $\mathfrak g_R$ is
\[
 [X,Y]_\beta=[X,Y]_0+\sum_{s=1}^r\kappa_s\Gamma^{(s)}(X,Y),
\]
the sum of the single-parameter deformations: $\Gamma^{(s)}$ is the coefficient $\Gamma_\beta$
for $\beta^{(s)}$. It is trivial: the even $R_r$-linear map
$U=\Theta_*T^{(1)}\cdots T^{(r)}$ satisfies $[UX,UY]_0=U([X,Y]_\beta)$ for all
$X,Y\in\mathfrak g_R$. Here $T^{(s)}$ is the trivialization of Theorem~\ref{thm:triv-single} for $\beta^{(s)}$,
$\Theta=(\mathbb I+\Omega J)^{1/2}$ with $\Omega_{uv}=-\sum_i\gamma_{u,i}\gamma_{v,i^*}$, and
$\Theta_*$ is the map that the substitution $B\mapsto\Theta B$ of the bosonic generators induces
on the lifts
(Theorems~\ref{thm:untwist-multi} and~\ref{thm:multi}, Proposition~\ref{prop:factor}).
\end{mainthm}

The trivialization is not additive, and this is the new feature for $r\geq2$. The product
$T^{(1)}\cdots T^{(r)}$ carries the deformation onto an even deformation, in which the bosonic
form $J$ acquires the central term $-\Omega$ of second order. The map $\Theta_*$, induced by the nonsymplectic matrix $\Theta$,
removes it. When $\Omega\neq0$, neither the product nor the truncation of $U$ to
degree $1$ is a trivialization (Proposition~\ref{prop:not-additive}). The statement holds for every
$N$ and $n$, including $D(n+1,n)$, whose Killing form vanishes. The Fock module of~\cite{Aoi26} extends to this setting, with an explicit vacuum (Remark~\ref{rem:fock}).

For $r=1$, Theorem~\ref{thm:A} extends~\cite{Aoi26} from $N=1$ to every $N$, and we establish
the following further results.

\begin{mainthm}\label{thm:B}
Let $r=1$. Then $U=\mathrm{id}+\kappa\tilde h_R$, where $\tilde h$ is the change of lifts
(Proposition~\ref{prop:lift-change}), and $\Gamma_\beta=\delta\tilde h$ for the coboundary
$\delta$ of \eqref{eq:coboundary}
(Theorems~\ref{thm:recovery} and~\ref{thm:triv-single}). The odd primitives of $\Gamma_\beta$
are exactly the maps $\tilde h-\ad Z$ with $Z\in P\otimes\mathfrak g_{\bar1}$, each for a unique
$Z$. Among them, $f_\beta=\tilde h-\ad Z_\beta$ is the only one vanishing on
$\mathfrak{sp}(2n)$, and its values are explicit (Lemma~\ref{lem:inner},
Theorems~\ref{thm:fbeta} and~\ref{thm:all-primitives}, Corollary~\ref{cor:normalization}).
For $N=1$, $f_\beta$ is the primitive of~\cite[Theorem~4.1]{Aoi26}.
\end{mainthm}

Two features are new for $N\geq2$. The odd lifts acquire a correction in
$\mathfrak{so}(N)$ when the source is untwisted (Proposition~\ref{prop:lift-change}).
Moreover, the canonical primitive $f_\beta$ does not vanish on $\mathfrak{so}(N)$. By
Corollary~\ref{cor:normalization} no primitive vanishes on the whole even part, so
the property ``$f_\beta$ vanishes on the even part'' of the case $N=1$ is special to
that case. These features need no separate treatment: the argument is uniform in $N$ and $n$ and is proved by index rules, with no case
analysis in the rank.

The structure of this paper is as follows. Section~\ref{sec:setting} fixes the index sets, the coefficients, the source algebra $A_\beta$ and the lifts.
Section~\ref{sec:untwist-base} untwists $A_\beta$ into the undeformed algebra $A_0$, where the
brackets of the lifts are computed and their span is identified with $\osp(N|2n)$. Section~\ref{sec:single} treats one odd
parameter: the change of lifts, the recovery and uniqueness of the coefficient, its explicit
form, the triviality and the primitives. Section~\ref{sec:multi} treats multiple odd
parameters, and Section~\ref{sec:special-remarks} the specializations and remarks.
Appendix~\ref{sec:formalization} records what is machine-checked.

For every $N$ the untwisting, the recovery and its uniqueness (in the algebras presented
by \eqref{eq:source}), the coboundary identity, the primitive $f_\beta$, the classification of the
primitives, and the identification with $\osp(N|2n)$ are machine-checked in Lean~4.
For every number $r$ of odd parameters, so are Theorems~\ref{thm:untwist-multi}
and~\ref{thm:multi} and Propositions~\ref{prop:lift-change-multi}, \ref{prop:factor} and~\ref{prop:not-additive}. What is not
formalized, including the Fock module of Remark~\ref{rem:fock}, is listed in
Appendix~\ref{sec:formalization}.

\section{The source algebra and the lifts}\label{sec:setting}

Fix $N\geq1$ and $n\geq1$, and write $m=\lfloor N/2\rfloor$.

\paragraph{Indices and forms.} Bosonic indices run over $I=\{1,\ldots,2n\}$. Let
\[
 J=\operatorname{diag}(J_2,\ldots,J_2),\qquad
 J_2=\begin{pmatrix}0&1\\-1&0\end{pmatrix}.
\]
Fermionic indices run over
\[
 K_N=\{1^+,1^-,\ldots,m^+,m^-\}\ (N=2m),\qquad
 K_N=\{0,1^+,1^-,\ldots,m^+,m^-\}\ (N=2m+1).
\]
For $k=1,\ldots,m$ put $(k^\pm)^*=k^\mp$, and put $0^*=0$. Let $g$ be the symmetric form on
$K_N$ with $g_{ii^*}=1$ for every $i$ and all other entries $0$; that is,
$g_{ij}=\delta_{i^*j}$ for all $i,j\in K_N$. Order $K_N$ as $0,1^+,1^-,2^+,\ldots$ (with $0$ omitted when $N$ is even).

\paragraph{Coefficients.} Put
\[
 P=\mathbb Q[\beta_{u,i}:u\in I,\ i\in K_N],\qquad
 R=P\oplus\kappa P,\qquad \kappa^2=0.
\]
All $\beta_{u,i}$ are even and $\kappa$ is odd, so $R$ is supercommutative.
The elements of $R$ are the scalars; they are written on the left. In an associative superalgebra (such as $A_\beta$
and $A_0$ below, or $\operatorname{End}(M)$ in Proposition~\ref{prop:osp}), the bracket without
subscript is the supercommutator: for homogeneous elements, with $p(\cdot)$ the
parity, $[X,Y]=XY-(-1)^{p(X)p(Y)}YX$. Braces denote this bracket when both arguments
are odd. Central means supercentral: $z$ is central if $[z,x]=0$ for all $x$; for even $z$ this
is the usual notion. Brackets on $\mathfrak g_R$ carry a subscript: $[\cdot,\cdot]_0$
(Section~\ref{sec:untwist-base}) and $[\cdot,\cdot]_\beta$ (Section~\ref{sec:recovery}).

\paragraph{The source.} Let $A_\beta$ be the unital associative $R$-superalgebra
generated by even $b_u$ ($u\in I$) and odd $a_i$ ($i\in K_N$), subject to
\begin{equation}\label{eq:source}
 [b_u,b_v]=J_{uv},\qquad \{a_i,a_j\}=g_{ij},\qquad [b_u,a_i]=\beta_{u,i}\kappa .
\end{equation}
The scalars are central; in particular $a_i\kappa=-\kappa a_i$. Thus
$a_0^2=\tfrac12$ when $N$ is odd, and $(a_{k^\pm})^2=0$. For $N=1$ this is the source
of~\cite{Aoi26}.

Let $A_0$ be the same algebra with $\beta=0$. The subscript $0$ marks $\beta=0$, here and in $\iota_0$ and $[\cdot,\cdot]_0$ below. We write $B_u$ for its bosonic
generators, so that $[B_u,a_i]=0$ in $A_0$. It is $R\otimes_{\mathbb Q}A_{\mathbb Q}$ with
$A_{\mathbb Q}=W_n\otimes_{\mathbb Q}\mathrm{Cl}(N)$, where $W_n$ is the rational Weyl algebra
on the $B_u$ and $\mathrm{Cl}(N)$ the Clifford algebra of $g$ on the $a_i$. (In~\cite{Aoi26}, $A_0$ and $A_{\mathbb Q}$ are written $A_B$ and $A_0$.) Tensor
products use graded signs. Thus $A_\beta$ is a deformation of $A_0$ over $R$: at
$\beta=0$ the two algebras coincide. The deformation of $\osp(N|2n)$ studied here is the one that $A_\beta$ induces
on the bracket of the lifts, $[\cdot,\cdot]_\beta$ (Section~\ref{sec:single}).

\paragraph{Lifts.} Define the symmetrized elements
\begin{equation}\label{eq:lifts}
 \widehat L_{uv}=\tfrac14(b_ub_v+b_vb_u),\qquad
 \widehat S_{ij}=\tfrac12(a_ia_j-a_ja_i),\qquad
 \widehat F_{u,i}=\tfrac14(a_ib_u+b_ua_i),
\end{equation}
for $u,v\in I$ and $i,j\in K_N$, with $\widehat L_{vu}=\widehat L_{uv}$, $\widehat S_{ji}=-\widehat S_{ij}$ and
$\widehat S_{ii}=0$. Write $L^0_{uv}$, $S^0_{ij}$, $F^0_{u,i}$ for the same
expressions in $A_0$. The symbols $L_{uv}$
($u\leq v$), $S_{ij}$ ($i<j$) and $F_{u,i}$ are a basis of a free $\mathbb Q$-module
$\mathfrak g$; as with the lifts, $L_{vu}=L_{uv}$, $S_{ji}=-S_{ij}$ and $S_{ii}=0$. We put $\mathfrak g_P=P\otimes\mathfrak g$ and
$\mathfrak g_R=R\otimes\mathfrak g$. Let $\iota_\beta$ and $\iota_0$ be the
$R$-linear maps from $\mathfrak g_R$ to $A_\beta$ and $A_0$ sending each symbol to its
lift. The $L$ and $S$ are even and the $F$ odd, so that
\[
 \dim\mathfrak g=\bigl(n(2n+1)+\tfrac{N(N-1)}{2}\ \big|\ 2nN\bigr),
\]
the superdimension of $\osp(N|2n)$~\cite[\S2.27]{FSS00}. We write $\mathfrak g_{\bar0}$ for the span of the $L$
and $S$, and $\mathfrak g_{\bar1}$ for the span of the $F$.

\begin{remark}[Even-central specialization]\label{rem:central}
The constants in \eqref{eq:source} impose $K=1$, as in~\cite[Remark~2.1]{Aoi26}, to which we
refer for this specialization. The odd parameter $\kappa$ is retained.
\end{remark}

\section{The untwisting and the identification with \texorpdfstring{$\osp(N|2n)$}{osp(N|2n)}}\label{sec:untwist-base}

We first show that the deformation $A_\beta$ of $A_0$ defined in Section~\ref{sec:setting} can be
removed by a change of generators (Lemma~\ref{lem:untwist}). We then work in $A_0$: we compute the
brackets of the lifts (Lemma~\ref{lem:base}), show that the lifts are independent
(Proposition~\ref{prop:independence}), and identify their span with $\osp(N|2n)$
(Proposition~\ref{prop:osp}). The change of generators does not carry lifts to lifts. Comparing
the lifts in $A_\beta$ and $A_0$ in Section~\ref{sec:single} produces the deformation of the bracket.

In $A_\beta$ the generators $b_u$ and $a_i$ fail to commute only by the odd central scalars
$\beta_{u,i}\kappa$. Shifting each $b_u$ by an even element linear in the $a_i$ absorbs these scalars.

\begin{lemma}[Untwisting]\label{lem:untwist}
Put $c_u=\kappa\sum_{i\in K_N}\beta_{u,i^*}a_i\in A_0$, and write $c_u$ also for the same expression in
$A_\beta$. There are mutually inverse even
isomorphisms of $R$-superalgebras $\Phi\colon A_\beta\to A_0$ and $\Psi\colon A_0\to A_\beta$,
both fixing the $a_i$ and $R$, with
\[
 \Phi(b_u)=B_u+c_u,\qquad \Psi(B_u)=b_u-c_u .
\]
\end{lemma}

\begin{proof}
Each $c_u$ is even, and $c_uc_v=0$ because $\kappa$ is odd and central with
$\kappa^2=0$. In $A_0$,
\[
 [B_u+c_u,a_j]=\kappa\sum_i\beta_{u,i^*}\{a_i,a_j\}
 =\kappa\sum_i\beta_{u,i^*}\,g_{ij}=\beta_{u,j}\kappa,
\]
and $[B_u+c_u,B_v+c_v]=J_{uv}+[B_u,c_v]+[c_u,B_v]=J_{uv}$, since $[B_u,a_i]=0$.
The fermionic relations and $R$ are unchanged. So the images of the generators satisfy
the relations \eqref{eq:source}, and $\Phi$ is well defined. Symmetrically $\Psi$ is well defined, using
$[b_u,c_v]=\kappa\sum_i\beta_{v,i^*}[b_u,a_i]=0$, which holds because $\kappa^2=0$.
Both maps fix every $c_u$, so $\Psi\Phi$ and $\Phi\Psi$ are the identity on
generators.
\end{proof}

For $N=1$, Lemma~\ref{lem:untwist} reduces to~\cite[Lemma~2.2]{Aoi26}.

We now work in $A_0$, where the bosons and the fermions supercommute. There the
supercommutator of two lifts is a linear combination of lifts.

\begin{lemma}[The base bracket]\label{lem:base}
In $A_0$ the lifts satisfy, for all indices,
\begin{align}
 [L^0_{uv},L^0_{wz}]&=\tfrac12\bigl(J_{vw}L^0_{uz}+J_{uw}L^0_{vz}+J_{vz}L^0_{uw}
   +J_{uz}L^0_{vw}\bigr),\label{eq:base-ll}\\
 [S^0_{ij},S^0_{kl}]&=g_{jk}S^0_{il}-g_{ik}S^0_{jl}-g_{jl}S^0_{ik}+g_{il}S^0_{jk},
   \label{eq:base-ss}\\
 [L^0_{uv},S^0_{ij}]&=0,\label{eq:base-ls}\\
 [L^0_{uv},F^0_{w,i}]&=\tfrac12\bigl(J_{vw}F^0_{u,i}+J_{uw}F^0_{v,i}\bigr),
   \label{eq:base-lf}\\
 [S^0_{jk},F^0_{w,i}]&=g_{ki}F^0_{w,j}-g_{ji}F^0_{w,k},\label{eq:base-sf}\\
 \{F^0_{u,i},F^0_{v,j}\}&=\tfrac14J_{uv}S^0_{ij}+\tfrac12g_{ij}L^0_{uv}.
   \label{eq:base-ff}
\end{align}
The remaining brackets follow by super-skew symmetry. No scalar term occurs.
\end{lemma}

\begin{proof}
In $A_0$ bosons and fermions supercommute, so
$F^0_{u,i}=\tfrac12a_iB_u$, $S^0_{ij}=a_ia_j-\tfrac12g_{ij}$ and
$L^0_{uv}=\tfrac12B_uB_v-\tfrac14J_{uv}$. The first line is~\cite[Lemma~2.3]{Aoi26}. The lines with $S$ follow from
$[a_ia_j,a_k]=g_{jk}a_i-g_{ik}a_j$. The lines with one $F$ follow from
$[L^0_{uv},B_w]=\tfrac12(J_{vw}B_u+J_{uw}B_v)$ and the previous rule. Finally
\[
 \{a_iB_u,a_jB_v\}=a_ia_jB_uB_v+a_ja_iB_vB_u .
\]
Substituting $a_ia_j=S^0_{ij}+\tfrac12g_{ij}$, $a_ja_i=-S^0_{ij}+\tfrac12g_{ij}$,
$B_uB_v=2L^0_{uv}+\tfrac12J_{uv}$ and $B_vB_u=2L^0_{uv}-\tfrac12J_{uv}$ gives
$J_{uv}S^0_{ij}+2g_{ij}L^0_{uv}$: the two constants $\tfrac14g_{ij}J_{uv}$ cancel.
Multiplying by $\tfrac14$ gives \eqref{eq:base-ff}.
\end{proof}

Write $[\cdot,\cdot]_0$ for the bracket on $\mathfrak g$ given by these rules, with $L^0$, $S^0$,
$F^0$ replaced by $L$, $S$, $F$, and extended to
$\mathfrak g_R$ by the sign rule of a supercommutator with central scalars:
\[
 [\lambda X,\mu Y]_0=(-1)^{p(X)p(\mu)}\lambda\mu[X,Y]_0\qquad(\lambda,\mu\in R,\ X,Y\ \text{homogeneous}).
\]
It is the
supercommutator of the lifts in $A_0$, and $\iota_0$ is injective
(Proposition~\ref{prop:independence}). So $[\cdot,\cdot]_0$ satisfies super-skew symmetry and the
super-Jacobi identity, and $\iota_0$ is a homomorphism: $\iota_0[X,Y]_0=[\iota_0X,\iota_0Y]$.
We write $\ad X$ for $Y\mapsto[X,Y]_0$. That no scalar term occurs means that the
supercommutator of two lifts is a combination of lifts, with no multiple of $1$. The
bracket $[\cdot,\cdot]_0$ is therefore defined on $\mathfrak g$ itself, not on a central
extension $\mathfrak g\oplus R\cdot1$.

\begin{proposition}[Independence]\label{prop:independence}
The map $\iota_0\colon\mathfrak g_R\to A_0$ is injective.
\end{proposition}

\begin{proof}
The ordered monomials in the $B_u$ form a basis of $W_n$~\cite[proof of Lemma~2.3]{Aoi26}, and the
ordered products of distinct $a_i$ form a basis of $\mathrm{Cl}(N)$. So the ordered monomials, fermions first,
form an $R$-basis of $A_0$. The lifts have distinct leading monomials: $B_uB_v$ for
$L^0_{uv}$, $a_ia_j$ for $S^0_{ij}$, and $a_iB_u$ for $F^0_{u,i}$. Hence they are
independent over $R$, so $\iota_0$ is injective.
\end{proof}

It remains to identify $(\mathfrak g,[\cdot,\cdot]_0)$. The rules
\eqref{eq:base-ll}--\eqref{eq:base-ff} are those of the natural action on a superspace with
an orthosymplectic form, which we now make explicit.

\begin{proposition}[The orthosymplectic identification]\label{prop:osp}
Let $M$ have even basis $e_i$ ($i\in K_N$) and odd basis $e_u$ ($u\in I$). In this
proposition and its proof, $i,j,k$ run over $K_N$ and $u,v,w$ over $I$. Let
$\mathcal B$ be the supersymmetric form with $\mathcal B(e_i,e_k)=-g_{ik}$,
$\mathcal B(e_u,e_w)=J_{uw}$ and $\mathcal B(e_i,e_u)=0$. Define $\rho\colon\mathfrak g\to
\operatorname{End}(M)$ by
\[
 \rho(L_{uv})\colon e_w\mapsto\tfrac12(J_{vw}e_u+J_{uw}e_v),\qquad
 \rho(S_{ij})\colon e_k\mapsto g_{jk}e_i-g_{ik}e_j,
\]
\[
 \rho(F_{u,i})\colon e_k\mapsto\tfrac12g_{ik}e_u,\quad e_w\mapsto\tfrac12J_{uw}e_i,
\]
with $\rho(L_{uv})$ zero on the $e_i$ and $\rho(S_{ij})$ zero on the $e_u$. Then $\rho$
is an isomorphism of $(\mathfrak g,[\cdot,\cdot]_0)$ onto $\osp(M,\mathcal B)$, the span of the
homogeneous $T\in\operatorname{End}(M)$ with
\[
 \mathcal B(Tx,y)+(-1)^{p(T)p(x)}\mathcal B(x,Ty)=0\qquad\text{for all homogeneous }x,y\in M.
\]
In the basis $e_{1^+},\ldots,e_{m^+},e_{1^-},\ldots,e_{m^-}$ ($e_0$ last if $N$ is odd),
then $e_2,e_4,\ldots,e_{2n},e_1,e_3,\ldots,e_{2n-1}$, the Gram matrix of $-\mathcal B$ is
the matrix $G$ of~\cite[\S2.27]{FSS00}.
\end{proposition}

\begin{proof}
\emph{The image lies in $\osp(M,\mathcal B)$.} This is a direct check on basis vectors, using that $g$
is symmetric and $J$ antisymmetric. For instance, $T=\rho(F_{u,i})$ is odd, and for
$x=e_k$, $y=e_w$ the condition reads
\[
 \mathcal B(Te_k,e_w)+\mathcal B(e_k,Te_w)=\tfrac12g_{ik}J_{uw}-\tfrac12J_{uw}g_{ki}=0 .
\]

\emph{Brackets.} Direct compositions give exactly \eqref{eq:base-ll}--\eqref{eq:base-ff}
under $\rho$.

\emph{Bijectivity.} The even part of $\osp(M,\mathcal B)$ preserves $M_{\bar0}$ and
$M_{\bar1}$, the spans of the $e_i$ and of the $e_u$. So it is the direct sum of the
endomorphisms $T$ of $M_{\bar0}$ with $\mathcal B(Tx,y)+\mathcal B(x,Ty)=0$ for $x,y\in M_{\bar0}$,
which form $\mathfrak{so}(N)$, and of the endomorphisms $T$ of $M_{\bar1}$ with the same condition
for $x,y\in M_{\bar1}$, which form $\mathfrak{sp}(2n)$. The $\rho(S_{ij})$ ($i<j$)
are a basis of the first summand. The $\rho(L_{uv})$ ($u\leq v$) are a basis of the
second: $T\mapsto((x,y)\mapsto J(Tx,y))$ identifies it with the symmetric forms, and $2$
is invertible. An odd element of $\osp(M,\mathcal B)$ is determined by its values on the
$e_i$, because the condition and nondegeneracy fix its values on $M_{\bar1}$. So the odd
part has dimension $2nN$, and the $\rho(F_{u,i})$ are a basis of it. Hence $\rho$ maps a
basis to a basis.

\emph{Normal form.} In the stated order, $-\mathcal B$ has even block
$\begin{pmatrix}0&\mathbb I_m\\ \mathbb I_m&0\end{pmatrix}$, with the entry $1$ added
when $N$ is odd, and odd block $\begin{pmatrix}0&\mathbb I_n\\-\mathbb I_n&0\end{pmatrix}$.
This is the matrix $G$ of~\cite[\S2.27]{FSS00}, over $\mathbb Q$ and with no change of sign of
basis vectors.
\end{proof}

We therefore identify $\mathfrak g$ with $\osp(N|2n)$. For $N=1$ this is the
identification of~\cite[Section~2]{Aoi26}. Under it the span of the $L_{uv}$ is $\mathfrak{sp}(2n)$ and that
of the $S_{ij}$ is $\mathfrak{so}(N)$, and we use these names for the two subspaces of $\mathfrak g$.

\section{The deformed bracket for one odd parameter}\label{sec:single}

This section treats one odd parameter, $r=1$, and extends~\cite{Aoi26} from $N=1$ to every
$N$. Section~\ref{sec:recovery} compares the lifts before and after the untwisting and
recovers from them the deformed bracket $[\cdot,\cdot]_\beta$. Its coefficient $\Gamma_\beta$ is the
coboundary $\delta\tilde h$ of the change of lifts, is determined by the lifts alone, and is computed
explicitly. Section~\ref{sec:triviality} derives the triviality from the primitive $\tilde h$.
Subtracting an inner derivation normalizes $\tilde h$ to $f_\beta$, which for $N=1$ is the primitive
of~\cite{Aoi26}. We then determine every odd primitive: $f_\beta$ is the unique one vanishing on
$\mathfrak{sp}(2n)$.

\subsection{The lifts and the coefficient}\label{sec:recovery}

The untwisting $\Phi$ of Lemma~\ref{lem:untwist} does not carry lifts to lifts; the difference
is linear in $\kappa$.

\begin{proposition}\label{prop:lift-change}
For every basis element $X$ of $\mathfrak g$, $\Phi(\iota_\beta X)=\iota_0(X+\kappa\tilde h(X))$, where $\tilde h$
is the odd $P$-linear map
\[
 \tilde h(L_{uv})=\sum_i\bigl(\beta_{u,i^*}F_{v,i}+\beta_{v,i^*}F_{u,i}\bigr),\qquad
 \tilde h(S_{ij})=0,\qquad
 \tilde h(F_{u,i})=\tfrac12\sum_j\beta_{u,j^*}S_{ji}.
\]
\end{proposition}

\begin{proof}
The map $\Phi$ fixes the $a_i$, so $\tilde h(S_{ij})=0$. For $L_{uv}$: $c_uc_v=0$, and $c_v$ commutes
with $B_u$. So
\[
 \Phi(\widehat L_{uv})-L^0_{uv}=\tfrac12(c_vB_u+c_uB_v)
 =\tfrac12\kappa\sum_i\bigl(\beta_{v,i^*}a_iB_u+\beta_{u,i^*}a_iB_v\bigr),
\]
and $\tfrac12a_iB_u=F^0_{u,i}$. For $F_{u,i}$: since
$a_i\kappa=-\kappa a_i$,
\[
 \Phi(\widehat F_{u,i})-F^0_{u,i}=\tfrac14(a_ic_u+c_ua_i)
 =\tfrac14\kappa\sum_j\beta_{u,j^*}(a_ja_i-a_ia_j),
\]
and $a_ja_i-a_ia_j=2S^0_{ji}$. This is a commutator of two fermions, so it has
no multiple of $1$, and $j=i$ contributes $S^0_{ii}=0$.
\end{proof}

Note that the correction $\tilde h(F_{u,i})$ lies in $P\otimes\mathfrak{so}(N)$. It vanishes for $N=1$,
where $K_1=\{0\}$ and $S_{00}=0$, and then $\tilde h$ is the map $h$ of~\cite[Section~3]{Aoi26}.

We express the effect of $\tilde h$ on the brackets through the coboundary of $[\cdot,\cdot]_0$. Let $\delta$ be the coboundary of an odd $P$-linear map $f\colon\mathfrak g_P\to
\mathfrak g_P$:
\begin{equation}\label{eq:coboundary}
 \delta f(X,Y)=[fX,Y]_0+(-1)^{p(X)}[X,fY]_0-f([X,Y]_0).
\end{equation}
Here an odd $P$-linear map $f$ is extended to $\mathfrak g_R$ by $f_R(\lambda X)=(-1)^{p(\lambda)}\lambda f(X)$.
An odd $P$-linear map $f$ with $\delta f=\Gamma$ is called a \emph{primitive} of a $P$-bilinear map
$\Gamma\colon\mathfrak g_P\times\mathfrak g_P\to\mathfrak g_P$. The brackets of the lifts
in $A_\beta$ follow by carrying those in $A_0$ back along $\Psi$.

\begin{theorem}[Recovery]\label{thm:recovery}
For all basis elements $X,Y$ of $\mathfrak g$,
\[
 [\iota_\beta X,\iota_\beta Y]=\iota_\beta\bigl([X,Y]_0+\kappa\,\delta\tilde h(X,Y)\bigr).
\]
In particular no $1$ or $\kappa\cdot1$ term occurs.
\end{theorem}

\begin{proof}
Two facts about $\kappa$, which is odd and central, are used: $[\kappa Y_1,Y_2]=\kappa[Y_1,Y_2]$
and $[Y_1,\kappa Y_2]=(-1)^{p(Y_1)}\kappa[Y_1,Y_2]$. By Proposition~\ref{prop:lift-change} and
Lemma~\ref{lem:base}, with $\kappa^2=0$ removing the product of the two corrections,
\[
 \Phi[\iota_\beta X,\iota_\beta Y]=[\iota_0(X+\kappa\tilde hX),\iota_0(Y+\kappa\tilde hY)]
 =\iota_0[X,Y]_0+\kappa\,\iota_0\bigl([\tilde hX,Y]_0+(-1)^{p(X)}[X,\tilde hY]_0\bigr).
\]
Proposition~\ref{prop:lift-change} gives $\Psi(\iota_0Z)=\iota_\beta(Z-\kappa\tilde hZ)$,
and $\Psi$ fixes $\kappa$. Applying $\Psi$ gives
\[
 [\iota_\beta X,\iota_\beta Y]=\iota_\beta[X,Y]_0+\kappa\,\iota_\beta\bigl([\tilde hX,Y]_0
 +(-1)^{p(X)}[X,\tilde hY]_0-\tilde h[X,Y]_0\bigr).
\]
This is the argument of~\cite[Section~3]{Aoi26}, with the map $\mathrm{id}+\kappa\tilde h_R$.
\end{proof}

\paragraph{Definition of the deformed bracket.} The map $\iota_\beta$ is injective: $\Phi\iota_\beta=\iota_0(\mathrm{id}+\kappa\tilde h_R)$
by Proposition~\ref{prop:lift-change}, the map $\mathrm{id}+\kappa\tilde h_R$ has inverse
$\mathrm{id}-\kappa\tilde h_R$, and $\iota_0$ is injective (Proposition~\ref{prop:independence}).
With Theorem~\ref{thm:recovery}, this defines the bracket $[\cdot,\cdot]_\beta$ on $\mathfrak g_R=R\otimes\mathfrak g$ (Section~\ref{sec:setting}), the
supercommutator of the lifts in $A_\beta$ carried to $\mathfrak g_R$:
\[
 \iota_\beta\bigl([X,Y]_\beta\bigr)=[\iota_\beta X,\iota_\beta Y]\qquad(X,Y\in\mathfrak g_R).
\]
Since $\iota_\beta$ is $R$-linear and the scalars are central, $[\cdot,\cdot]_\beta$ is $R$-bilinear with the
sign rule of $[\cdot,\cdot]_0$. For basis elements $X,Y$ of $\mathfrak g$,
Theorem~\ref{thm:recovery} and the injectivity of $\iota_\beta$ give
$[X,Y]_\beta=[X,Y]_0+\kappa\Gamma_\beta(X,Y)$ with
$\Gamma_\beta(X,Y)\in\mathfrak g_P$. This defines the \emph{coefficient} $\Gamma_\beta$, a $P$-bilinear
map $\mathfrak g_P\times\mathfrak g_P\to\mathfrak g_P$, and Theorem~\ref{thm:recovery} says
$\Gamma_\beta=\delta\tilde h$.

Thus $[\cdot,\cdot]_\beta$ is a deformation of $[\cdot,\cdot]_0$ on $\mathfrak g_R$ itself, with no central
term. In~\cite{BS17} the central elements $K$ and $\kappa$ appear at the level of the generators, but
by Theorem~\ref{thm:recovery} no term in $1$ or $\kappa\cdot1$ occurs at the level of the lifts.
The reason is the symmetrization \eqref{eq:lifts}: a scalar can only arise from reordering two
generators, and these scalars cancel in the supercommutator of symmetrized elements (as in
the proof of Lemma~\ref{lem:base} and in Proposition~\ref{prop:lift-change}).
The coefficient is determined by the lifts alone:

\begin{theorem}[Uniqueness]\label{thm:unique}
Let $X,Y$ be basis elements of $\mathfrak g$. If $[\iota_\beta X,\iota_\beta Y]=\iota_\beta([X,Y]_0)+\kappa\,\iota_\beta(\Gamma')$ with
$\Gamma'\in\mathfrak g_P$, then $\Gamma'=\Gamma_\beta(X,Y)$.
\end{theorem}

\begin{proof}
Subtracting Theorem~\ref{thm:recovery} gives $\kappa\,\iota_\beta(\Gamma'-\Gamma_\beta(X,Y))=0$.
By Proposition~\ref{prop:lift-change} and $\kappa^2=0$, applying $\Phi$ gives
$\kappa\,\iota_0(\Gamma'-\Gamma_\beta(X,Y))=0$. In $A_0=P\otimes A_{\mathbb Q}\oplus\kappa P\otimes A_{\mathbb Q}$,
multiplication by $\kappa$ is injective on $P\otimes A_{\mathbb Q}$. So
$\iota_0(\Gamma'-\Gamma_\beta(X,Y))=0$, and Proposition~\ref{prop:independence} gives
$\Gamma'=\Gamma_\beta(X,Y)$.
\end{proof}

Evaluating $\delta\tilde h$ on basis elements gives $\Gamma_\beta$ explicitly.

\begin{proposition}[The explicit coefficient]\label{prop:gamma}
For all indices,
\begin{align}
 \Gamma_\beta(L_{uv},L_{wz})&=0,\qquad \Gamma_\beta(S_{ij},S_{kl})=0,\label{eq:gamma-ll}\\
 \Gamma_\beta(L_{uv},S_{ij})&=\beta_{u,i}F_{v,j}-\beta_{u,j}F_{v,i}
   +\beta_{v,i}F_{u,j}-\beta_{v,j}F_{u,i},\label{eq:gamma-ls}\\
 \Gamma_\beta(L_{uv},F_{w,k})&=\tfrac12\bigl(\beta_{u,k}L_{vw}+\beta_{v,k}L_{uw}\bigr),
   \label{eq:gamma-lf}\\
 \Gamma_\beta(S_{ij},F_{w,k})&=\tfrac12\bigl(\beta_{w,j}S_{ik}-\beta_{w,i}S_{jk}\bigr),
   \label{eq:gamma-sf}\\
 \Gamma_\beta(F_{u,i},F_{v,j})&=-\tfrac12\bigl(\beta_{u,j}F_{v,i}+\beta_{v,i}F_{u,j}\bigr),
   \label{eq:gamma-ff}
\end{align}
and $\Gamma_\beta(Y,X)=-(-1)^{p(X)p(Y)}\Gamma_\beta(X,Y)$.
\end{proposition}

\begin{proof}
Expand $\delta\tilde h$ with Lemma~\ref{lem:base} and
Proposition~\ref{prop:lift-change}, using $g_{ij}=\delta_{i^*j}$. For example, for
$(L_{uv},F_{w,k})$ the terms $[\tilde hL_{uv},F_{w,k}]_0$ contribute
$\tfrac12(\beta_{u,k}L_{vw}+\beta_{v,k}L_{uw})$ through the $L$-part of
\eqref{eq:base-ff}, and their $S$-parts cancel against
$[L_{uv},\tilde hF_{w,k}]_0-\tilde h[L_{uv},F_{w,k}]_0$. The other kinds are similar.
\end{proof}

For $N=1$ only \eqref{eq:gamma-ll}, \eqref{eq:gamma-lf} and \eqref{eq:gamma-ff}
remain, and they are those of~\cite[Proposition~3.1]{Aoi26}. The coefficient $\Gamma_\beta$ itself is
nonzero, since \eqref{eq:gamma-lf} is nonzero.

\subsection{Triviality and the primitives}\label{sec:triviality}

The map $\tilde h$ of Proposition~\ref{prop:lift-change} is a primitive of $\Gamma_\beta$
(Theorem~\ref{thm:recovery}), and this gives the triviality.

\begin{theorem}[Triviality for one odd parameter]\label{thm:triv-single}
The maps $T=\mathrm{id}+\kappa\tilde h_R$ and $T^{-1}=\mathrm{id}-\kappa\tilde h_R$ are even,
$R$-linear, mutually inverse, and satisfy
\[
 [TX,TY]_0=T([X,Y]_\beta)
\]
for all $X,Y\in\mathfrak g_R$, where $[\cdot,\cdot]_\beta$ is the bracket of Section~\ref{sec:recovery}.
\end{theorem}

\begin{proof}
The argument of~\cite[Theorem~4.1]{Aoi26} applies with $\tilde h$ in place of $f_\beta$. For
homogeneous basis elements, $\kappa^2=0$ and $[X,\kappa Y]_0=(-1)^{p(X)}\kappa[X,Y]_0$ reduce the
identity to $\Gamma_\beta=\delta\tilde h$, which is Theorem~\ref{thm:recovery}. Since $T$ is even and
$R$-linear, both sides are $R$-bilinear with the sign rule of $[\cdot,\cdot]_0$, so the identity extends to $\mathfrak g_R$. Since $(\kappa\tilde h_R)^2=0$, the inverse is
$\mathrm{id}-\kappa\tilde h_R$.
\end{proof}

The primitive $\tilde h$ does not vanish on $\mathfrak{sp}(2n)$. Subtracting an inner derivation
corrects this. Define $\check\beta_{t,i}\in P$ by $\sum_t\check\beta_{t,i}J_{tu}=\beta_{u,i}$, that is,
\begin{equation}\label{eq:v}
 \check\beta_{2k-1,i}=\beta_{2k,i},\qquad \check\beta_{2k,i}=-\beta_{2k-1,i}\qquad(k=1,\ldots,n,\ i\in K_N),
\end{equation}
and put
\begin{equation}\label{eq:primitive}
 Z_\beta=2\sum_{t,i}\check\beta_{t,i^*}F_{t,i}\in P\otimes\mathfrak g_{\bar1},\qquad
 f_\beta=\tilde h-\ad Z_\beta .
\end{equation}

\begin{lemma}\label{lem:inner}
For every $Z\in P\otimes\mathfrak g_{\bar1}$, $\delta(\ad Z)=0$.
\end{lemma}

\begin{proof}
The ring $P$ is even and central, so $[\cdot,\cdot]_0$ is $P$-bilinear with no signs and the super-Jacobi
identity extends from $\mathfrak g$ to $\mathfrak g_P$. Since $Z$ is odd with coefficients in $P$,
$\ad Z$ is odd and $P$-linear. For homogeneous $X,Y$,
super-Jacobi with $p(Z)=1$ reads
$[Z,[X,Y]_0]_0=[[Z,X]_0,Y]_0+(-1)^{p(X)}[X,[Z,Y]_0]_0$, which is
$\delta(\ad Z)(X,Y)=0$.
\end{proof}

\begin{theorem}\label{thm:fbeta}
We have $\Gamma_\beta=\delta f_\beta$, and for all indices
\[
 f_\beta(L_{uv})=0,\qquad
 f_\beta(F_{w,j})=-\sum_t\check\beta_{t,j}L_{tw},\qquad
 f_\beta(S_{ij})=2\sum_t\bigl(\check\beta_{t,j}F_{t,i}-\check\beta_{t,i}F_{t,j}\bigr).
\]
So $f_\beta$ vanishes on $\mathfrak{sp}(2n)$, maps the odd part into $\mathfrak{sp}(2n)$
and maps $\mathfrak{so}(N)$ into the odd part. For $N=1$, $f_\beta$ is the primitive of~\cite[Theorem~4.1]{Aoi26}.
\end{theorem}

\begin{proof}
By Lemma~\ref{lem:inner}, $\delta f_\beta=\delta\tilde h=\Gamma_\beta$. The values follow from
Lemma~\ref{lem:base}, super-skew symmetry and $\sum_t\check\beta_{t,i}J_{tu}=\beta_{u,i}$. For $L$:
\[
 [Z_\beta,L_{uv}]_0=-\sum_{t,i}\check\beta_{t,i^*}\bigl(J_{vt}F_{u,i}+J_{ut}F_{v,i}\bigr)
 =\sum_i\bigl(\beta_{v,i^*}F_{u,i}+\beta_{u,i^*}F_{v,i}\bigr)=\tilde h(L_{uv}).
\]
For $F$: $[F_{t,i},F_{w,j}]_0=\tfrac14J_{tw}S_{ij}+\tfrac12g_{ij}L_{tw}$, so
\[
 [Z_\beta,F_{w,j}]_0=\tfrac12\sum_i\beta_{w,i^*}S_{ij}+\sum_t\check\beta_{t,j}L_{tw},
\]
whose first sum is $\tilde h(F_{w,j})$. For $S$: $\tilde h(S_{ij})=0$, and
$[S_{ij},F_{t,l}]_0=g_{jl}F_{t,i}-g_{il}F_{t,j}$ gives
$f_\beta(S_{ij})=[S_{ij},Z_\beta]_0=2\sum_t(\check\beta_{t,j}F_{t,i}-\check\beta_{t,i}F_{t,j})$. For $N=1$
there is no $S$, and, writing $F_u=F_{u,0}$ and $\beta_u=\beta_{u,0}$ as in~\cite{Aoi26},
$f_\beta(F_u)=\sum_{k=1}^n(\beta_{2k-1}L_{(2k)u}-\beta_{2k}L_{(2k-1)u})$, which is~\cite[(4.2)]{Aoi26}.
\end{proof}

The $\mathfrak{so}(N)$-part of $\tilde h(F)$ is thus canceled exactly by that of
$[Z_\beta,F]_0$. For $N=1$ the relation $f_\beta=\tilde h-\ad Z_\beta$ is shown
after~\cite[Theorem~4.1]{Aoi26}, where it is written $f_\beta=h-\ad(2F(v))$, with $h=\tilde h$ (as noted after Proposition~\ref{prop:lift-change}) and $v_t=\check\beta_{t,0}$. Here it is the definition of $f_\beta$, \eqref{eq:primitive}.

\begin{corollary}\label{cor:tbeta}
The map $T_\beta=\mathrm{id}+\kappa(f_\beta)_R$ has inverse $\mathrm{id}-\kappa(f_\beta)_R$ and
satisfies $[T_\beta X,T_\beta Y]_0=T_\beta([X,Y]_\beta)$ for all $X,Y\in\mathfrak g_R$.
\end{corollary}

\begin{proof}
The proof is the same as for Theorem~\ref{thm:triv-single}, with $\Gamma_\beta=\delta f_\beta$.
\end{proof}

We now have two odd primitives of $\Gamma_\beta$, $\tilde h$ and $f_\beta$ (Theorems~\ref{thm:recovery} and~\ref{thm:fbeta}). Primitives are determined only up to odd cocycles, so a trivialization
involves a choice. It is therefore natural to ask which odd primitives $\Gamma_\beta$ has. For $N=1$, where $\mathfrak{sp}(2n)$ is the even part, the
primitive of~\cite{Aoi26} vanishes there. Theorem~\ref{thm:all-primitives} answers the question,
and with it, vanishing on $\mathfrak{sp}(2n)$ determines $f_\beta$ for every $N$, while vanishing on
the whole even part is impossible for $N\geq2$ (Corollary~\ref{cor:normalization}).

\begin{theorem}\label{thm:all-primitives}
Let $f$ be an odd $P$-linear map $\mathfrak g_P\to\mathfrak g_P$ with $\delta f=\Gamma_\beta$.
Then $f=\tilde h-\ad Z$ for a unique $Z\in P\otimes\mathfrak g_{\bar1}$. Equivalently, every odd
$P$-linear map $D\colon\mathfrak g_P\to\mathfrak g_P$ with $\delta D=0$, that is, every odd $P$-linear
derivation of $\mathfrak g_P$, is $\ad Z$ for a unique $Z\in P\otimes\mathfrak g_{\bar1}$.
\end{theorem}

With $\mathbb C$ in place of $P$, the second statement is the classical fact that every derivation of
$\osp(N|2n)$ over $\mathbb C$ is inner~\cite[Chapter~III, \S1.2, Proposition~3]{Sch79}. For
$\osp(N|2n)\otimes A$, with $A$ unital, associative and supercommutative over an algebraically closed
field of characteristic zero, see~\cite[Corollary~4.13(2)]{Duf02}.

\begin{proof}
\emph{Reduction to $\mathbb Q$.} The map $D=f-\tilde h$ is odd and satisfies $\delta D=0$, that
is, $D([X,Y]_0)=[DX,Y]_0+(-1)^{p(X)}[X,DY]_0$. Expanding its values in the monomial
basis of $P$ writes $D=\sum_\mu\mu D_\mu$, a finite sum of odd $\mathbb Q$-linear maps.
The structure constants are rational, so each $D_\mu$ satisfies the same identity. It
suffices to show that every such $D$ is $\ad Z'$ for a unique $Z'\in\mathfrak g_{\bar1}$.

\emph{On $\mathfrak{sp}(2n)$.} Let $\mathfrak{sp}(2n)$ act on $\mathfrak g_{\bar1}$ by
$L\cdot y=[L,y]_0$. For $L,L'\in\mathfrak{sp}(2n)$ the identity gives
$D([L,L']_0)=L\cdot D(L')-L'\cdot D(L)$, so $D|_{\mathfrak{sp}(2n)}$ is a $1$-cocycle of
$\mathfrak{sp}(2n)\cong\mathfrak{sp}(2n,\mathbb Q)$ (Proposition~\ref{prop:osp}) with values in a finite-dimensional module.
By Whitehead's lemma~\cite[Chapter~III, \S7, Lemma~3]{Jac}, which holds over any field of
characteristic zero, there is $y\in\mathfrak g_{\bar1}$ with $D(L)=L\cdot y$ for all
$L$. (Jacobson writes modules on the right, $f([ab])=f(a)b-f(b)a$, and concludes
$f(a)=da$. With $L\cdot x=-xL$ his condition is the one above, and his conclusion reads
$D(L)=L\cdot y$ with $y=-d$.) Put $Z'=-y$. Since $y$ is odd and $L$ even,
$L\cdot y=[L,y]_0=-[y,L]_0=[Z',L]_0$, so $D|_{\mathfrak{sp}(2n)}=\ad Z'|_{\mathfrak{sp}(2n)}$.

\emph{Equivariance.} The map $D'=D-\ad Z'$ is odd, satisfies $\delta D'=0$ by
Lemma~\ref{lem:inner}, and vanishes on $\mathfrak{sp}(2n)$. So
$D'([L,X]_0)=[L,D'X]_0$ for $L\in\mathfrak{sp}(2n)$: $D'$ is $\mathfrak{sp}(2n)$-equivariant.

\emph{Vanishing.} As $\mathfrak{sp}(2n)$-modules, $\mathfrak{sp}(2n)$ is the adjoint module, simple of
dimension $n(2n+1)$; $\mathfrak{so}(N)$ is trivial by \eqref{eq:base-ls}; and
$\mathfrak g_{\bar1}\cong(\mathbb Q^{2n})^{\oplus N}$, where for each fixed $i$ the span of the
$F_{u,i}$ is the natural module $\mathbb Q^{2n}$, simple of dimension $2n$, by
\eqref{eq:base-lf}. The map $D'|_{\mathfrak{so}(N)}\colon\mathfrak{so}(N)\to\mathfrak g_{\bar1}$ takes values in the
invariants of $(\mathbb Q^{2n})^{\oplus N}$, which are $0$. The map
$D'|_{\mathfrak g_{\bar1}}\colon(\mathbb Q^{2n})^{\oplus N}\to\mathfrak{sp}(2n)\oplus\mathfrak{so}(N)$ is $0$ by
Schur's lemma, since $2n\neq n(2n+1)$ and $\mathbb Q^{2n}$ is not trivial. The natural module is simple over every field containing $\mathbb Q$, by a direct check with $L_{uu}$
and $J$ invertible. The algebra $\mathfrak{sp}(2n)\otimes\mathbb C=\mathfrak{sp}(2n,\mathbb C)$ is simple~\cite{Jac}, so the
adjoint module is simple over $\mathbb Q$.
Hence $D'=0$.

\emph{Uniqueness.} If $\ad Z'=\ad Z''$, then $Z'-Z''$ commutes with $\mathfrak{sp}(2n)$, so it
lies in the invariants of $\mathfrak g_{\bar1}$, which are $0$. This proves the last statement, and the
first follows with $D=f-\tilde h$.
\end{proof}

\begin{corollary}[Normalization]\label{cor:normalization}
The map $f_\beta$ is the unique odd primitive of $\Gamma_\beta$ vanishing on $\mathfrak{sp}(2n)$.
For $N\geq2$, no odd primitive vanishes on $\mathfrak{sp}(2n)\oplus\mathfrak{so}(N)$.
\end{corollary}

\begin{proof}
By Theorem~\ref{thm:all-primitives} a primitive is $\tilde h-\ad Z$. It vanishes on
$\mathfrak{sp}(2n)$ exactly when $\ad(Z-Z_\beta)|_{\mathfrak{sp}(2n)}=0$, that is, when $Z=Z_\beta$. So a
primitive vanishing on $\mathfrak{sp}(2n)\oplus\mathfrak{so}(N)$ would be $f_\beta$.
For $N\geq2$ take $i<j$ in $K_N$: the coefficient of $F_{1,i}$ in $f_\beta(S_{ij})$ is
$2\check\beta_{1,j}=2\beta_{2,j}\neq0$. Hence $f_\beta(S_{ij})\neq0$, and no odd primitive vanishes on
$\mathfrak{sp}(2n)\oplus\mathfrak{so}(N)$.
\end{proof}

\begin{corollary}[Complex parameters]\label{cor:complex}
For any values $\beta_{u,i}\in\mathbb C$, the specialized bracket on
$(\mathbb C\oplus\kappa\mathbb C)\otimes\mathfrak g$ is isomorphic to the undeformed bracket
via $\mathrm{id}+\kappa f_\beta$.
\end{corollary}

\begin{proof}
The proof is the same as for~\cite[Corollary~4.2]{Aoi26}: the identities of Theorem~\ref{thm:fbeta} and
Corollary~\ref{cor:tbeta} are polynomial in $\beta$ with rational coefficients, and the inverse
is explicit.
\end{proof}

\section{The deformed bracket for multiple odd parameters}\label{sec:multi}

With one odd parameter, the bracket $[\cdot,\cdot]_0+\kappa\Gamma$ satisfies the Jacobi identity for
every odd $2$-cocycle $\Gamma$, such as $\Gamma_\beta=\delta\tilde h$. Its parts of degree $0$ and $1$ in
$\kappa$ are the Jacobi identity of $[\cdot,\cdot]_0$ and the cocycle condition, and its part of
degree $2$ vanishes as $\kappa^2=0$. So no second-order obstruction arises. On the other hand, deformations can be defined over
any supercommutative base~\cite[\S4.1]{Lei23}, such as the Grassmann algebra
$\Lambda(\kappa_1,\ldots,\kappa_r)$, the functions on the superpoint $0|r$~\cite[\S2.1]{Lei23}. Over it
the products $\kappa_s\kappa_{s'}$ with $s\neq s'$ are even and nonzero, and second-order effects appear.
In fact, they appear only in the trivialization: the deformed bracket stays linear in the $\kappa_s$
(Theorem~\ref{thm:multi}), while its trivialization is not additive and at second order involves an
even central perturbation of the bosonic form (Propositions~\ref{prop:factor} and~\ref{prop:not-additive}, Remark~\ref{rem:even-generated}).

We now allow $r\geq1$ odd parameters; the new phenomena need $r\geq2$. Put
\[
 P_r=\mathbb Q[\beta^{(s)}_{u,i}:1\leq s\leq r,\ u\in I,\ i\in K_N],\qquad
 R_r=P_r\otimes\Lambda(\kappa_1,\ldots,\kappa_r),
\]
with the $\kappa_s$ odd. The $\kappa$-degree of a product of $\kappa_s$ is the number of its factors.
Every element of $R_r$ is a sum of such products with coefficients in $P_r$. It has $\kappa$-degree at most $r$, since a product of more than $r$ of the $\kappa_s$ repeats one of them. Let $A_\beta$ be generated
over $R_r$ as in
Section~\ref{sec:setting}, with the mixed relation of \eqref{eq:source} replaced by
\begin{equation}\label{eq:source-multi}
 [b_u,a_i]=\gamma_{u,i},\qquad \gamma_{u,i}=\sum_{s=1}^r\beta^{(s)}_{u,i}\kappa_s .
\end{equation}
The $\gamma_{u,i}$ are odd and central, and $\gamma\gamma'=-\gamma'\gamma$. We keep
$A_0$, the lifts, $\mathfrak g_R$, $\iota_\beta$ and $\iota_0$, now over $R_r$. In particular $\mathfrak g_R=R_r\otimes\mathfrak g$. For $r=1$ this is the
setting of Section~\ref{sec:setting}, with $\kappa=\kappa_1$ and $\beta=\beta^{(1)}$, so that $P_1=P$ and $R_1=R$.

\subsection{Recovery and triviality}\label{sec:multi-triv}

As in the case $r=1$, we first untwist. Put $c_u=\sum_i\gamma_{u,i^*}a_i\in A_0$, which is even. As in the proof of
Lemma~\ref{lem:untwist}, $[c_u,a_j]=\gamma_{u,j}$ and $[B_w,c_u]=0$. But $c_uc_v$ need
not vanish once $r\geq2$:
\begin{equation}\label{eq:Omega}
 \Omega_{uv}:=[c_u,c_v]=-\sum_i\gamma_{u,i}\gamma_{v,i^*}.
\end{equation}
Indeed $[c_u,c_v]=-\sum_{i,j}\gamma_{u,i^*}\gamma_{v,j^*}\{a_i,a_j\}$ and
$g_{ij}=\delta_{i^*j}$, after reindexing. The matrix $\Omega$ is antisymmetric, its entries are even,
central and of $\kappa$-degree $2$, and $\Omega=0$ when $r=1$.

\begin{theorem}\label{thm:untwist-multi}
Let $\Theta=(\mathbb I+\Omega J)^{1/2}=\sum_{k\geq0}\binom{1/2}{k}(\Omega J)^k$, a finite sum ending at
$k=\lfloor r/2\rfloor$, since the entries of $(\Omega J)^k$ have $\kappa$-degree $2k$. There are
mutually inverse even isomorphisms of $R_r$-superalgebras $\Phi_r\colon A_\beta\to A_0$ and
$\Psi_r\colon A_0\to A_\beta$, both fixing the $a_i$ and $R_r$, with
\[
 \Phi_r(b_u)=\sum_w\Theta_{uw}B_w+c_u,\qquad
 \Psi_r(B_u)=\sum_w(\Theta^{-1})_{uw}(b_w-c_w).
\] In $A_\beta$, $c_w$ denotes the
same expression $\sum_i\gamma_{w,i^*}a_i$.
\end{theorem}

\begin{proof}
Since $\Omega^{\mathsf T}=-\Omega$ and $J^{\mathsf T}=-J$, one has
$(\Omega J)^{\mathsf T}=J\Omega$, so $J(\Omega J)^{\mathsf T}=J^2\Omega=-\Omega$ and $(\Omega J)J=\Omega J^2=-\Omega$. Hence
$J\varphi((\Omega J)^{\mathsf T})=\varphi(\Omega J)J$ for every polynomial $\varphi$, and
$\Theta J\Theta^{\mathsf T}=\Theta^2J=(\mathbb I+\Omega J)J=J-\Omega$. In $A_0$ the $B$ commute with the $c$, so
\[
 [\Phi_rb_u,\Phi_rb_v]=(\Theta J\Theta^{\mathsf T})_{uv}+\Omega_{uv}=J_{uv},\qquad
 [\Phi_rb_u,a_j]=[c_u,a_j]=\gamma_{u,j}.
\]
As in the proof of Lemma~\ref{lem:untwist}, the fermionic relations and $R_r$ are unchanged.
Hence the images of the generators satisfy \eqref{eq:source-multi}, and $\Phi_r$ is well defined. In $A_\beta$ put $b'_w=b_w-c_w$. There
$[b_w,c_v]=\sum_i\gamma_{v,i^*}\gamma_{w,i}=\Omega_{wv}$ and $[c_w,c_v]=\Omega_{wv}$, so
$[b'_w,b'_v]=J_{wv}-\Omega_{wv}$ and $[b'_w,a_j]=\gamma_{w,j}-\gamma_{w,j}=0$. Since
$J-\Omega=\Theta J\Theta^{\mathsf T}$, the elements $\Psi_rB_u$ satisfy the Weyl relations and commute
with the $a_j$, so $\Psi_r$ is well defined. Finally $\Psi_r\Phi_r(b_u)=b'_u+c_u=b_u$,
and $\Phi_r\Psi_r(B_u)=B_u$ in the same way, since $\Phi_r(c_w)=c_w$.
\end{proof}

In particular $\Theta=\mathbb I+\tfrac12\Omega J$ for $r\leq3$. For $r=1$, $\Theta=\mathbb I$ and
Theorem~\ref{thm:untwist-multi} reduces to Lemma~\ref{lem:untwist}.

As in the case $r=1$, the untwisting does not carry lifts to lifts.

\begin{proposition}\label{prop:lift-change-multi}
For every basis element $X$ of $\mathfrak g$, $\Phi_r(\iota_\beta X)=\iota_0(UX)$, where $U$ is the even
$R_r$-linear map with $U(S_{ij})=S_{ij}$,
\begin{align*}
 U(L_{uv})&=\sum_{w,z}\Theta_{uw}\Theta_{vz}L_{wz}
  +\sum_{w,i}\bigl(\Theta_{uw}\gamma_{v,i^*}+\Theta_{vw}\gamma_{u,i^*}\bigr)F_{w,i}
  -\tfrac12\sum_{i,j}\gamma_{u,i^*}\gamma_{v,j^*}S_{ij},\\
 U(F_{u,i})&=\sum_w\Theta_{uw}F_{w,i}+\tfrac12\sum_j\gamma_{u,j^*}S_{ji}.
\end{align*}
The map $U-\mathrm{id}$ raises the $\kappa$-degree, so it is nilpotent and $U$ is invertible.
\end{proposition}

\begin{proof}
As in the proof of Proposition~\ref{prop:lift-change}, the map $\Phi_r$ fixes the $a_i$, so
$U(S_{ij})=S_{ij}$. Substitute $\Phi_rb_u=(\Theta B)_u+c_u$ into \eqref{eq:lifts}. The part in $(\Theta B)(\Theta B)$ gives
the first sum. The cross terms are $F$-lifts, since the $B$ commute with the $a$. The
remaining terms are $\tfrac14(c_uc_v+c_vc_u)=-\tfrac12\sum_{i,j}\gamma_{u,i^*}\gamma_{v,j^*}S^0_{ij}$
in $\widehat L_{uv}$, and $\tfrac14(a_ic_u+c_ua_i)=\tfrac12\sum_j\gamma_{u,j^*}S^0_{ji}$
in $\widehat F_{u,i}$. Both are commutators of fermions, so no multiple of $1$ occurs.
\end{proof}

For $r=1$, $\Theta=\mathbb I$ and $U=\mathrm{id}+\kappa\tilde h_R$, as in Proposition~\ref{prop:lift-change}.

\begin{theorem}[Triviality for $r$ odd parameters]\label{thm:multi}
For all basis elements $X,Y$ of $\mathfrak g$,
\[
 [\iota_\beta X,\iota_\beta Y]=\iota_\beta\bigl(U^{-1}[UX,UY]_0\bigr),
\]
so the supercommutator of two lifts in $A_\beta$ is an $R_r$-linear combination of lifts, with no
scalar term, and $\iota_\beta$ is injective. Thus
$[\cdot,\cdot]_\beta$ is defined on $\mathfrak g_R$ as in the case $r=1$ (Section~\ref{sec:recovery}), and the even $R_r$-linear map $U$ satisfies
$[UX,UY]_0=U([X,Y]_\beta)$ for all $X,Y\in\mathfrak g_R$.
Moreover
\begin{equation}\label{eq:linear}
 [X,Y]_\beta=[X,Y]_0+\sum_{s=1}^r\kappa_s\,\Gamma^{(s)}(X,Y),
\end{equation}
where $\Gamma^{(s)}$ is the coefficient $\Gamma_\beta$ of Proposition~\ref{prop:gamma}
with $\beta^{(s)}$ in place of $\beta$.
\end{theorem}

\begin{proof}
By Proposition~\ref{prop:lift-change-multi}, $\Phi_r\iota_\beta=\iota_0U$, so
$\Psi_r\iota_0=\iota_\beta U^{-1}$. Since $\iota_0$ is a homomorphism for
$[\cdot,\cdot]_0$, applying $\Psi_r$ to
$\Phi_r[\iota_\beta X,\iota_\beta Y]=[\iota_0UX,\iota_0UY]=\iota_0[UX,UY]_0$ gives
the first identity. The map $\iota_0$ is injective over $R_r$ by the proof of
Proposition~\ref{prop:independence}, hence so is $\iota_\beta=\Psi_r\iota_0U$. By the injectivity of $\iota_\beta$,
the first identity says $[X,Y]_\beta=U^{-1}[UX,UY]_0$ for basis elements $X,Y$, that is,
$[UX,UY]_0=U([X,Y]_\beta)$. As in the proof of Theorem~\ref{thm:triv-single}, the identity extends
to $\mathfrak g_R$, since $U$ is even and $R_r$-linear.

For \eqref{eq:linear}, compute $[\iota_\beta X,\iota_\beta Y]$ in $A_\beta$ directly.
Both lifts are quadratic in the generators. The super-Leibniz rule
$[x,yz]=[x,y]z+(-1)^{p(x)p(y)}y[x,z]$ of the supercommutator, for homogeneous $x,y,z$, applied in each
argument, shows that the bracket is a sum of terms $\pm x_1[x_2,x_3]x_4$, each with exactly one bracket of two generators. That
bracket is $J_{uv}$, $g_{ij}$ or $\pm\gamma_{u,i}$, a central scalar of $\kappa$-degree
at most $1$. A product $x_1x_4$ of two generators is
$\tfrac12\bigl(x_1x_4+(-1)^{p(x_1)p(x_4)}x_4x_1\bigr)$, a rational multiple of a lift~\eqref{eq:lifts},
plus $\tfrac12[x_1,x_4]$, again such a scalar. By the first part the scalar terms cancel, since
$1$ and the lifts are independent in $A_\beta$ (by $\Phi_r$ and the ordered monomial
basis of the proof of Proposition~\ref{prop:independence}). The coefficients of the
lifts are then the entries of $[X,Y]_\beta$. The scalars $\tfrac12[x_1,x_4]$ contribute only to the scalar part. So each coefficient of a lift
is a rational combination of single brackets $[x_2,x_3]$. We carry the $\gamma_{u,i}$
through as formal odd central symbols. Then each entry is
$c^0+\sum_{u,i}\gamma_{u,i}C_{u,i}$, with rational $c^0$ and $C_{u,i}$ independent of $r$
and $\beta$. The signs do not depend on $r$, since each $\gamma_{u,i}$ is odd. At $\beta=0$, where $A_\beta=A_0$, Lemma~\ref{lem:base} gives
$c^0=[X,Y]_0$. At $r=1$, where $\gamma=\beta\kappa$, comparing the coefficients of $\kappa$ with
Theorem~\ref{thm:recovery} gives
$\Gamma_\beta=\sum_{u,i}\beta_{u,i}C_{u,i}$. In general, substituting
$\gamma_{u,i}=\sum_s\kappa_s\beta^{(s)}_{u,i}$ gives the entry $[X,Y]_0+\sum_s\kappa_s\sum_{u,i}\beta^{(s)}_{u,i}C_{u,i}$.
The identity $\Gamma_\beta=\sum_{u,i}\beta_{u,i}C_{u,i}$ is polynomial in the indeterminates $\beta_{u,i}$. So
$\sum_{u,i}\beta^{(s)}_{u,i}C_{u,i}$ is the entry of $\Gamma^{(s)}(X,Y)$, which gives \eqref{eq:linear}.
\end{proof}

For $r=1$, $U=\mathrm{id}+\kappa\tilde h_R$ is the map $T$ of Theorem~\ref{thm:triv-single}, and
Theorem~\ref{thm:multi} reduces to Theorems~\ref{thm:recovery} and~\ref{thm:triv-single}.

So the deformed bracket is the sum of the single-parameter deformations of
Section~\ref{sec:single}, and each is trivial. The trivialization, however, is not
the sum of their trivializations. It is their product, followed by a change of the bosonic generators.

\subsection{Factorization through an even deformation}\label{sec:factor}

For a family $\chi=(\chi_{u,i})$ of odd central elements of $R_r$, let $V_\chi$ be the even
$R_r$-linear map with $V_\chi(S_{ij})=S_{ij}$,
\begin{align*}
 V_\chi(L_{uv})&=L_{uv}+\sum_i\bigl(\chi_{v,i^*}F_{u,i}+\chi_{u,i^*}F_{v,i}\bigr)
  -\tfrac12\sum_{i,j}\chi_{u,i^*}\chi_{v,j^*}S_{ij},\\
 V_\chi(F_{u,i})&=F_{u,i}+\tfrac12\sum_j\chi_{u,j^*}S_{ji},
\end{align*}
that is, the map $U$ of Proposition~\ref{prop:lift-change-multi} with $\Theta=\mathbb I$ and $\chi$ in
place of $\gamma$. Put $T^{(s)}=V_{\beta^{(s)}\kappa_s}$. Since $\kappa_s^2=0$, the term quadratic in $\chi$ vanishes, and
$T^{(s)}=\mathrm{id}+\kappa_s\tilde h^{(s)}_R$, where $\tilde h^{(s)}$ is the map $\tilde h$ of Proposition~\ref{prop:lift-change} for $\beta^{(s)}$. So $T^{(s)}$ is the trivialization of Theorem~\ref{thm:triv-single} for $\beta^{(s)}$.
Let $\Theta_*$ be the even $R_r$-linear map with $\Theta_*(S_{ij})=S_{ij}$,
$\Theta_*(L_{uv})=\sum_{w,z}\Theta_{uw}\Theta_{vz}L_{wz}$ and $\Theta_*(F_{u,i})=\sum_w\Theta_{uw}F_{w,i}$.
Let $[\cdot,\cdot]_{J-\Omega}$ be the bracket on $\mathfrak g_R$ given by
\eqref{eq:base-ll}--\eqref{eq:base-ff} with $J-\Omega$ in place of $J$, extended to
$\mathfrak g_R$ as $[\cdot,\cdot]_0$ is.

\begin{proposition}\label{prop:factor}
The $T^{(s)}$ commute with each other, and
\[
 U=\Theta_*T^{(1)}\cdots T^{(r)}.
\]
The product $T^{(1)}\cdots T^{(r)}$ intertwines $[\cdot,\cdot]_\beta$ with $[\cdot,\cdot]_{J-\Omega}$, and
$\Theta_*$ intertwines $[\cdot,\cdot]_{J-\Omega}$ with $[\cdot,\cdot]_0$. If $\Omega\neq0$, then $\Theta_*$ is not an
automorphism of $[\cdot,\cdot]_0$.
\end{proposition}

\begin{proof}
First, $V_\chi V_\psi=V_{\chi+\psi}$. On $F_{u,i}$ this is clear, since $V$ fixes the $S$.
On $L_{uv}$, apply $V_\chi$ to $V_\psi(L_{uv})$. The terms $\psi_{v,i^*}F_{u,i}$ contribute
$\tfrac12\sum_{i,j}\psi_{v,i^*}\chi_{u,j^*}S_{ji}=-\tfrac12\sum_{i,j}\chi_{u,i^*}\psi_{v,j^*}S_{ij}$,
since odd scalars anticommute. The terms $\psi_{u,i^*}F_{v,i}$ contribute the same sum with $u$ and $v$ exchanged. With the $S$-terms of
$V_\chi$ and $V_\psi$ this gives the $S$-term of $V_{\chi+\psi}$. As
$\chi+\psi=\psi+\chi$, the maps $V_\chi$ commute, and $T^{(1)}\cdots T^{(r)}=V_\gamma$.

In $A_\beta$ put $b'_u=b_u-c_u$. By the proof of Theorem~\ref{thm:untwist-multi},
$[b'_u,b'_v]=J_{uv}-\Omega_{uv}$, $[b'_u,a_j]=0$ and $\Phi_r(b'_u)=\sum_w\Theta_{uw}B_w$. Let $\iota'$ send
each symbol to its lift \eqref{eq:lifts} written with $b'$ in place of $b$. The relations of
Lemma~\ref{lem:base} are formal consequences of the brackets of the generators. Here these are
$[b'_u,b'_v]=J_{uv}-\Omega_{uv}$, $[b'_u,a_i]=0$ and $\{a_i,a_j\}=g_{ij}$, which are central scalars. So
$[\iota'X,\iota'Y]=\iota'[X,Y]_{J-\Omega}$. Substituting $\Phi_r(b')=\Theta B$ into the lifts gives
$\Phi_r\iota'=\iota_0\Theta_*$. The same formulas define $\Pi_*$ for every matrix $\Pi$ with even central entries, and
$\Pi_*(\Pi')_*=(\Pi'\Pi)_*$. Since $\Theta$ is $\mathbb I$ plus a nilpotent matrix, $\Theta_*$ is invertible, with inverse
$(\Theta^{-1})_*$. So
$\iota'$ is injective, since $\iota_0$ is. Substituting $b=b'+c$ as in the proof of Proposition~\ref{prop:lift-change-multi}, with
$b'$ in place of $\Theta B$, gives $\iota_\beta=\iota'V_\gamma$. Hence $\iota_0U=\Phi_r\iota_\beta=\iota_0\Theta_*V_\gamma$, so
$U=\Theta_*V_\gamma$. For the second intertwining, applying $\Phi_r$ to $[\iota'X,\iota'Y]=\iota'[X,Y]_{J-\Omega}$ gives
$[\iota_0\Theta_*X,\iota_0\Theta_*Y]=\iota_0\Theta_*[X,Y]_{J-\Omega}$. Since $\iota_0$ is an injective homomorphism,
$[\Theta_*X,\Theta_*Y]_0=\Theta_*[X,Y]_{J-\Omega}$. For the first,
$\iota_\beta=\iota'V_\gamma$ and Theorem~\ref{thm:multi} give
$\iota'V_\gamma[X,Y]_\beta=[\iota'V_\gamma X,\iota'V_\gamma Y]=\iota'[V_\gamma X,V_\gamma Y]_{J-\Omega}$. Since $\iota'$ is injective,
$V_\gamma[X,Y]_\beta=[V_\gamma X,V_\gamma Y]_{J-\Omega}$.

Suppose $\Theta_*$ is an automorphism of $[\cdot,\cdot]_0$. Then, for all $X,Y$,
$\Theta_*[X,Y]_{J-\Omega}=[\Theta_*X,\Theta_*Y]_0=\Theta_*[X,Y]_0$, and $[X,Y]_{J-\Omega}=[X,Y]_0$ since $\Theta_*$ is
injective. But \eqref{eq:base-ll} with $v=u$ and $z=w$ gives $[L_{uu},L_{ww}]_0=2J_{uw}L_{uw}$ and
$[L_{uu},L_{ww}]_{J-\Omega}=2(J_{uw}-\Omega_{uw})L_{uw}$, so the difference is $-2\Omega_{uw}L_{uw}$. If $\Omega_{uw}\neq0$, then $u\neq w$, since
$\Omega$ is antisymmetric. So $L_{uw}$ is a basis element, and the difference is nonzero. This contradicts
$[X,Y]_{J-\Omega}=[X,Y]_0$, so $\Theta_*$ is not an automorphism of $[\cdot,\cdot]_0$.
\end{proof}

In particular, if $\Omega\neq0$, then $U(T^{(1)}\cdots T^{(r)})^{-1}=\Theta_*$ is not of the
form $\exp(\ad Y)$ for any even $Y\in\mathfrak g_R$ with $\ad Y$ nilpotent: such a map is an automorphism of $[\cdot,\cdot]_0$. For $r\leq3$, $\Theta=\mathbb I+\tfrac12\Omega J$, so the map $\Theta_*-\mathrm{id}$ has $\kappa$-degree $2$ and its square
vanishes. In that case, if $\Omega\neq0$, the map $\Theta_*-\mathrm{id}$ is not a derivation of $[\cdot,\cdot]_0$, since otherwise
$\Theta_*=\exp(\Theta_*-\mathrm{id})$ would be an automorphism.

\begin{proposition}[The trivialization is not additive]\label{prop:not-additive}
Suppose $\Omega\neq0$. Then neither $T^{(1)}\cdots T^{(r)}$ nor the truncation
$U_{\le1}=\mathrm{id}+\sum_s\kappa_s\tilde h^{(s)}_R$ of $U$ to degree $1$ intertwines
$[\cdot,\cdot]_\beta$ with $[\cdot,\cdot]_0$.
\end{proposition}

\begin{proof}
By Theorem~\ref{thm:multi}, $U$ intertwines $[\cdot,\cdot]_\beta$ with $[\cdot,\cdot]_0$. If
$T^{(1)}\cdots T^{(r)}$ did as well, then $\Theta_*=U(T^{(1)}\cdots T^{(r)})^{-1}$ would be an
automorphism of $[\cdot,\cdot]_0$, contrary to Proposition~\ref{prop:factor}. For $U_{\le1}$, we have
$r\geq2$, since $\Omega=0$ when $r=1$. Take $X=L_{11}$ and $Y=L_{22}$. By
\eqref{eq:linear} and \eqref{eq:gamma-ll}, $[X,Y]_\beta=[X,Y]_0=2L_{12}$, so $U_{\le1}([X,Y]_\beta)$ has
no $\kappa_1\kappa_2$-part. But $[\kappa_1Y_1,\kappa_2Y_2]_0=-\kappa_1\kappa_2[Y_1,Y_2]_0$ for odd $Y_1,Y_2$, and
\eqref{eq:base-ff} then gives the $L_{12}$-coefficient
\[
 -2\sum_i\bigl(\beta^{(1)}_{1,i^*}\beta^{(2)}_{2,i}-\beta^{(2)}_{1,i^*}\beta^{(1)}_{2,i}\bigr)
\]
of the $\kappa_1\kappa_2$-part of $[U_{\le1}X,U_{\le1}Y]_0$. It is twice the coefficient of $\kappa_1\kappa_2$ in $\Omega_{12}$,
and it is nonzero since the $\beta$ are indeterminate. Hence $[U_{\le1}X,U_{\le1}Y]_0\neq U_{\le1}([X,Y]_\beta)$.
\end{proof}

\begin{remark}[The odd parameters generate an even perturbation]\label{rem:even-generated}
Proposition~\ref{prop:factor} shows where the even perturbation enters. The product of the
single-parameter trivializations carries the odd deformation onto $[\cdot,\cdot]_{J-\Omega}$. At
second order this is the even central perturbation $-\Omega$ of the bosonic form $J$, in the
even nilpotent parameters $\kappa_s\kappa_{s'}$. The map $\Theta_*$ removes it, and $\Theta$ is not symplectic:
$\Theta J\Theta^{\mathsf T}=J-\Omega$ by the proof of Theorem~\ref{thm:untwist-multi}. In
$[\Phi_rb_u,\Phi_rb_v]$ the same cancellation appears as the commutator $[c_u,c_v]=\Omega_{uv}$
of the shifts $c_u$. The even direction of Remark~\ref{rem:even} is thus reached from the odd
one, although it is not treated here.
\end{remark}

\begin{remark}[The Fock module]\label{rem:fock}
In this remark $k,l$ run over $\{1,\ldots,n\}$.
Let $\mathcal F=R_r[x_1,\ldots,x_n]\otimes_{\mathbb Q}\mathrm{Cl}(N)$, with the $x_k$ even and $\mathrm{Cl}(N)$
acting on itself by left multiplication. Sending $B_{2k-1}\mapsto\partial_{x_k}$ and $B_{2k}\mapsto x_k$
represents $A_0$ on $\mathcal F$, and composing with $\Phi_r$ represents $A_\beta$, faithfully, as in the case $N=1$, $r=1$ of~\cite{Aoi26}. For $N=n=r=1$ this is the Fock space of~\cite[Section~4]{BS17}, where the mixed relation
is realized by an odd derivation and the annihilated element is $1$ (see~\cite{Aoi26}). Put
$\eta_{kl}=\Theta_{(2k-1)(2l-1)}$ and $\zeta_{kl}=\Theta_{(2k-1)(2l)}$. Since $\Theta$ is $\mathbb I$ plus a matrix with entries of
$\kappa$-degree at least $2$, so is $\eta$, and $\eta$ is invertible. Put $\Xi=\eta^{-1}\zeta$ and $\tilde c_k=\sum_l(\eta^{-1})_{kl}c_{2l-1}$, and let $\Xi^+$ and $\Xi^-$ be the
symmetric and antisymmetric parts of $\Xi$. The elements annihilated by every $b_{2k-1}$ are
exactly the $\varpi\,\sigma$ with $\sigma\in R_r\otimes\mathrm{Cl}(N)$, where
\[
 \varpi=\exp\Bigl(-\tfrac12\sum_{k,l}\Xi^+_{kl}x_kx_l-\sum_kx_k\tilde c_k\Bigr),
\]
a finite sum. Indeed, $\mathfrak b_k=\sum_l(\eta^{-1})_{kl}b_{2l-1}$ acts as
$\partial_{x_k}+\sum_l\Xi_{kl}x_l+\tilde c_k$, and the $\mathfrak b_k$ commute, since $J_{(2k-1)(2l-1)}=0$. This gives
$[\tilde c_k,\tilde c_l]=2\Xi^-_{kl}$, a central scalar. With $\xi=\sum_kx_k\tilde c_k$, one has
$\partial_{x_k}e^{-\xi}=-(\tilde c_k-\tfrac12[\xi,\tilde c_k])e^{-\xi}$. The Gaussian factor
$\exp(-\tfrac12\sum_{k,l}\Xi^+_{kl}x_kx_l)$ is central in $\mathcal F$, and $\partial_{x_k}$ multiplies it by $-(\Xi^+x)_k$, which absorbs the
symmetric part of $\Xi$. The $\tilde c_k$ act by left multiplication. The antisymmetric part $\Xi^-$ cancels
the term $\tfrac12[\xi,\tilde c_k]=-(\Xi^-x)_k$, which arises because the $\tilde c_k$ do not commute. So $\mathfrak b_k\varpi=0$. Conversely, $\varpi$ is invertible, and $\mathfrak b_k(\varpi\sigma')=\varpi\,\partial_{x_k}\sigma'$ for every
$\sigma'\in\mathcal F$. With $\Omega^{\mathrm{oo}}_{kl}=\Omega_{(2k-1)(2l-1)}$, $\Omega^{\mathrm{oe}}_{kl}=\Omega_{(2k-1)(2l)}$ and
$c_{\mathrm o}=(c_{2k-1})_k$, the table below lists the cases worked out. In it, $\theta$, $f$ and $\lambda$ are even central elements
of $R_r$.
\begin{center}\small
\begin{tabular}{llll}
case & $\Theta$ & $\Xi$ & $\varpi$\\\hline
$r=1$ & $\mathbb I$ & $0$ & $1-\sum_kx_kc_{2k-1}$\\
$r=2$ & $\mathbb I+\tfrac12\Omega J$ & $\tfrac12\Omega^{\mathrm{oo}}$ & $\exp(-\sum_kx_kc_{2k-1})$\\
$r=3$ & $\mathbb I+\tfrac12\Omega J$ & $\tfrac12\Omega^{\mathrm{oo}}$ & $\exp(-\sum_kx_k\tilde c_k)$, $\tilde c=(\mathbb I+\tfrac12\Omega^{\mathrm{oe}})c_{\mathrm o}$\\
$n=1$ & $\theta\,\mathbb I$ & $0$ & $\exp(-\theta^{-1}x_1c_1)$\\
$N=1$ & $\mathbb I+f\,\Omega J$ & $\lambda\,(\gamma_{2k-1,0}\gamma_{2l-1,0})_{kl}$ & $\exp(-\sum_kx_k\tilde c_k)$
\end{tabular}
\end{center}
For $r=2$, $\eta^{-1}=\mathbb I+\tfrac12\Omega^{\mathrm{oe}}$ and $\Omega^{\mathrm{oe}}c_{\mathrm o}$ has $\kappa$-degree $3$, so
$\tilde c=c_{\mathrm o}$. For $n=1$, $\Omega=\omega J$ and $\Theta=(1-\omega)^{1/2}\mathbb I$. For $N=1$, $(\Omega J)^2$ is an
even multiple of $\Omega J$, which gives the form of $\Theta$. The Gaussian factor
$\exp(-\tfrac12\sum_{k,l}\Xi^+_{kl}x_kx_l)$ is $1$ unless $N,n\geq2$ and $r\geq4$. Indeed, $\Xi^+=0$ in three cases. For $r\leq3$, $\Theta=\mathbb I+\tfrac12\Omega J$ and
$\Xi=\tfrac12(\Omega_{(2k-1)(2l-1)})_{kl}$, which is antisymmetric. For $n=1$, $\Theta$ is a multiple of
$\mathbb I$, so $\zeta=0$ and $\Xi=0$. For $N=1$, $\Omega_{uv}=-\gamma_{u,0}\gamma_{v,0}$ has rank~$1$, and
$\Xi_{kl}=\lambda\,\gamma_{2k-1,0}\gamma_{2l-1,0}$ with a single even $\lambda$ for all $k,l$, so $\Xi$ is
antisymmetric. For $N,n\geq2$ and $r\geq4$, $\Xi^+$ need not vanish. Its part of $\kappa$-degree $4$ is
$\tfrac18(\Omega^{\mathrm{oe}}\Omega^{\mathrm{oo}}-\Omega^{\mathrm{oo}}(\Omega^{\mathrm{oe}})^{\mathsf T})$. At $N=n=2$, $r=4$, with
$\gamma_{1,1^+}=\kappa_1$, $\gamma_{3,1^-}=\kappa_2$, $\gamma_{1,1^-}=\kappa_3$, $\gamma_{4,1^+}=\kappa_4$ and the other
$\gamma_{u,i}=0$, one has $\Omega_{13}=-\kappa_1\kappa_2$ and $\Omega_{14}=-\kappa_3\kappa_4$, and this gives
$\Xi^+_{11}=-\tfrac14\Omega_{13}\Omega_{14}=-\tfrac14\kappa_1\kappa_2\kappa_3\kappa_4$. The factor depends on the
choice of $\Theta$: replacing $\Theta$ by $\Theta\mathcal S$, with $\mathcal S_{(2k-1)(2l)}=q_{kl}$ for a symmetric
nilpotent even $q$ and $\mathcal S$ otherwise the identity, keeps $\Theta J\Theta^{\mathsf T}=J-\Omega$ and changes
$\Xi^+$ into $\Xi^++q$.
\end{remark}

\section{Specializations and remarks}\label{sec:special-remarks}

This section specializes the results to $N=1$, to the smallest case of each series with $N\geq2$, and to the
realizations of~\cite{FSS00}, and closes with two remarks.

\paragraph{$N=1$.} With $K_1=\{0\}$ there is no $S$, $\tilde h(F)=0$, and
Lemma~\ref{lem:untwist}, Lemma~\ref{lem:base}, Theorem~\ref{thm:recovery} and
Theorem~\ref{thm:fbeta} reduce to~\cite[Lemma~2.2, Lemma~2.3, Proposition~3.1,
Theorem~4.1]{Aoi26}.

\paragraph{The $\mathfrak{so}(N)$ correction.} For every $N\geq2$ the correction
$\tilde h(F)\in\mathfrak{so}(N)$ is nonzero: for $i=1^+$ the term $j=1^-$ of
Proposition~\ref{prop:lift-change} gives $\tfrac12\beta_{u,1^+}S_{1^-1^+}$. The smallest
cases of the series $C(n+1)$ and $B(m,n)$, $m\geq1$, are $C(2)=\osp(2|2)$ and $B(1,1)=\osp(3|2)$; that of
$D(m,n)$ is $D(2,1)$ (Remark~\ref{rem:d21}). In each, Theorem~\ref{thm:fbeta} gives
$f_\beta(\mathfrak{so}(N))\neq0$, as Corollary~\ref{cor:normalization} requires.

\paragraph{The oscillator realizations of~\cite{FSS00}.} Frappat, Sciarrino and Sorba
realize $B(m,n)$ with $m\geq1$, $C(n+1)$ and $D(m,n)$ by $n$ bosonic oscillator pairs
$b_k^\pm$ and $m$ fermionic pairs $a_i^\pm$, together with one real fermion $e$, $e^2=1$,
for $B(m,n)$~\cite[\S2.30, pp.~218--221]{FSS00}. In this paragraph, $b_k^\pm$, $a_i^\pm$, $e$, $H_i$ and
$E_{\pm\alpha_i}$ are their notation. Identify
\[
 b_k^-=B_{2k-1},\qquad b_k^+=B_{2k},\qquad a_i^+=a_{i^+},\qquad a_i^-=a_{i^-},\qquad
 e=\sqrt2\,a_0 ,
\]
with $k=1,\ldots,n$ and $i=1,\ldots,m$. Here $i$ and $j$ index the simple roots. Then each of their
simple generators $H_i$, $E_{\pm\alpha_i}$ in the distinguished basis
lies in $\iota_0(\mathfrak g)$, with no scalar remainder. In particular, their constants,
such as the $-1$ in $H_{n+m}=2a_m^+a_m^--1$ for $B(m,n)$, are exactly the constants of
the symmetrization \eqref{eq:lifts}. Moreover,
$H_i=[E_{\alpha_i},E_{-\alpha_i}]$. After the rescaling of~\cite[\S2.2]{FSS00}, which turns the $H_i$ into $h_i$, the matrix of $[h_i,E_{\alpha_j}]$ is the
distinguished Cartan matrix of~\cite[Tables~3.54, 3.56, 3.57]{FSS00}. The root vectors
they list lie in $\iota_0(\mathfrak g)$ and carry the roots they assign. Two misprints
in~\cite{FSS00} are corrected in this comparison: the highest root in Table~3.56 is
$\alpha_1+2\alpha_2+\cdots+2\alpha_n+\alpha_{n+1}=\varepsilon+\delta_1$, and the row of
Table~2.6 labeled $C(n)$ is that of $C(n+1)$. For $B(0,n)$, Frappat, Sciarrino and Sorba~\cite{FSS00} use bosonic
oscillators only, and the correspondence is that of~\cite[Section~2]{Aoi26}. This comparison
is an independent check of Proposition~\ref{prop:osp}; the proof of that proposition does not use it.

\begin{remark}[The even direction]\label{rem:even}
The construction of~\cite[Section~4]{BS17} has an even, $K$-valued part as well as the
odd, $\kappa$-valued part treated here: the homogeneous relations $[b_u,b_v]=KJ_{uv}$ and
$\{a_i,a_j\}=Kg_{ij}$, which this paper specializes to $K=1$ (Remark~\ref{rem:central}).
Deformations in $K$, with $K$ kept as an indeterminate, are not considered. The argument does,
however, give the following. For an even antisymmetric $\Omega$ with $\Omega J$ nilpotent,
the perturbed bosonic form $J-\Omega$ is equivalent to $J$: $\Theta J\Theta^{\mathsf T}=J-\Omega$
with $\Theta=(\mathbb I+\Omega J)^{1/2}$, by the proof of Theorem~\ref{thm:untwist-multi}. With
multiple odd parameters such a perturbation arises at second order from the odd one
(Remark~\ref{rem:even-generated}).
\end{remark}

\begin{remark}[$D(2,1)$]\label{rem:d21}
For $N=4$, $n=1$ the algebra is $D(2,1)=\osp(4|2)$, the member $\alpha=1$ of the family
$D(2,1;\alpha)$~\cite[\S2.20]{FSS00}. Its Killing form vanishes~\cite[\S2.23]{FSS00}, as does that of every $D(n+1,n)$. The results above hold for
these cases without change, since the argument does not use the Killing form of
$\osp(N|2n)$. In particular, Theorem~\ref{thm:all-primitives} uses Whitehead's lemma only for
$\mathfrak{sp}(2n)$. Since $\osp(4|2)$ is
the only $\osp(N|2n)$ of superdimension $(9|8)$, a member $D(2,1;\alpha)$ of the family lies in
the present setting exactly when it is isomorphic to $D(2,1)$, as it is for
$\alpha\in\{1,-2,-\tfrac12\}$~\cite[\S2.20]{FSS00}.
\end{remark}

\appendix

\section{Machine-checked formalization}\label{sec:formalization}

The statements below have been formalized and machine-checked in the Lean~4 proof
assistant for every $N\geq1$ and $n\geq1$, and those of
Section~\ref{sec:multi} for every $r\geq1$. The fermionic index is a finite
type with an involution $i\mapsto i^*$, instantiated at $K_N$ for odd and for even $N$.
The development is the repository~\cite{Dev}, which accompanies this paper. From the formalization~\cite{Repo} of the case $N=1$, it uses the Weyl-algebra and
graded-tensor-product infrastructure and, for the controls of item~9 below, the objects of that case.
The modules that provide them are copied into the repository with their statements and proofs
unchanged. The repository builds with no warnings. Nothing in the body of the paper depends on this appendix.

\subsection{What is proved}\label{subsec:formal-proved}

\begin{enumerate}
\item The base bracket $[\cdot,\cdot]_0$ of Lemma~\ref{lem:base} is super-skew-symmetric and
is compatible with the $\mathbb Z/2$-grading.
\item The identity $\Gamma_\beta=\delta\tilde h$ holds on $\mathfrak g_P$, with $\Gamma_\beta$ given by
Proposition~\ref{prop:gamma} and $\tilde h$ by Proposition~\ref{prop:lift-change}.
\item The map $T=\mathrm{id}+\kappa\tilde h_R$ intertwines $[\cdot,\cdot]_0+\kappa\Gamma_\beta$, which is
$[\cdot,\cdot]_\beta$ by item~4, with $[\cdot,\cdot]_0$ for arbitrary elements of $\mathfrak g_R$, with inverse $\mathrm{id}-\kappa\tilde h_R$
(Theorem~\ref{thm:triv-single}).
\item The algebras $A_\beta$ and $A_0$ are constructed as the superalgebras presented
by~\eqref{eq:source}. In them: Lemma~\ref{lem:untwist}, with $\Phi$ and $\Psi$ mutually inverse
even isomorphisms; Lemma~\ref{lem:base}; Proposition~\ref{prop:lift-change};
Theorem~\ref{thm:recovery}; and the definition of $[\cdot,\cdot]_\beta$ in
Section~\ref{sec:recovery}: the supercommutator of two lifts is the image under $\iota_\beta$ of a
unique element of $\mathfrak g_R$, with $[X,Y]_\beta=[X,Y]_0+\kappa\Gamma_\beta(X,Y)$; and
$\iota_\beta([X,Y]_\beta)=[\iota_\beta X,\iota_\beta Y]$ for $X,Y$ in the $P$-basis of $\mathfrak g_R$ formed
by the basis elements of $\mathfrak g$ and their products with $\kappa$, hence for all of $\mathfrak g_R$
by $P$-bilinearity.
\item Proposition~\ref{prop:independence}: $\iota_0$ is injective on $\mathfrak g_R$, and so is
$\iota_\beta$. The coefficient is unique (Theorem~\ref{thm:unique}).
\item Lemma~\ref{lem:inner}, the values of $f_\beta$ and $\Gamma_\beta=\delta f_\beta$
(Theorem~\ref{thm:fbeta}), Corollary~\ref{cor:tbeta}, Theorem~\ref{thm:all-primitives},
Corollary~\ref{cor:normalization}, and Corollary~\ref{cor:complex} for any values
$\beta_{u,i}\in\mathbb C$. For $N\geq2$, $\tilde h(F_{u,i})\neq0$ and $f_\beta(S_{ij})\neq0$ for
$i\neq j$ (Section~\ref{sec:special-remarks} and the proof of Corollary~\ref{cor:normalization}).
\item Proposition~\ref{prop:osp}: $\rho$ takes values in $\osp(M,\mathcal B)$, preserves brackets,
is injective, and has every homogeneous element of $\osp(M,\mathcal B)$ in its image, so that it is
bijective onto $\osp(M,\mathcal B)$. In the stated order, the Gram matrix of $-\mathcal B$ is the matrix $G$ of~\cite{FSS00}. Mathlib has no
Lie superalgebras, so $M$, $\mathcal B$ and $\osp(M,\mathcal B)$ are defined in the
development's own terms.
\item Multiple odd parameters. The coefficients $R_r$ are modeled with the Grassmann signs
of the products of the $\kappa_s$, and $A_\beta$ and $A_0$ are constructed as the
superalgebras presented by~\eqref{eq:source-multi}. In them:
\begin{enumerate}
\item Theorem~\ref{thm:untwist-multi}, with $\Phi_r$ and $\Psi_r$ mutually inverse even
isomorphisms. The matrix $\Theta=\sum_{k\leq\lfloor r/2\rfloor}\binom{1/2}{k}(\Omega J)^k$ is even,
has constant part $\mathbb I$ and satisfies $\Theta J\Theta^{\mathsf T}=J-\Omega$; for $r\leq3$ it
is $\mathbb I+\tfrac12\Omega J$;
\item Proposition~\ref{prop:lift-change-multi}: $\Phi_r(\iota_\beta X)=\iota_0(UX)$, and the map
$U$, defined for any even $\Theta$, is invertible when $\Theta$ has constant part $\mathbb I$;
\item Theorem~\ref{thm:multi}, with $\iota_\beta$ injective, and
$\iota_\beta([X,Y]_\beta)=[\iota_\beta X,\iota_\beta Y]$ for all homogeneous $X,Y\in\mathfrak g_R$. The identity
$[UX,UY]_0=U([X,Y]_\beta)$, for $[\cdot,\cdot]_\beta$ given by~\eqref{eq:linear}, holds for every
even $\Theta$ with $\Theta J\Theta^{\mathsf T}=J-\Omega$;
\item Proposition~\ref{prop:factor}: $V_\chi V_\psi=V_{\chi+\psi}$ for odd $\chi,\psi$, so
the $T^{(s)}$ commute and $T^{(1)}\cdots T^{(r)}=V_\gamma$; $\Pi_*(\Pi')_*=(\Pi'\Pi)_*$ for even
$\Pi,\Pi'$; $U=\Theta_*T^{(1)}\cdots T^{(r)}$; the two intertwinings through
$[\cdot,\cdot]_{J-\Omega}$; and, if $\Omega\neq0$, $\Theta_*$ is not an automorphism of
$[\cdot,\cdot]_0$;
\item the text after Proposition~\ref{prop:factor}: the super-Jacobi identity of
$[\cdot,\cdot]_0$ on basis elements of $\mathfrak g_R$, so that $\ad Y$ is a derivation for
every even $Y$, and $\exp(\ad Y)$ an automorphism when $\ad Y$ is nilpotent; hence, if
$\Omega\neq0$, $\Theta_*$ is not $\exp(\ad Y)$ for any such $Y$. For $r\leq3$ and
$\Omega\neq0$, $\Theta_*-\mathrm{id}$ is not a derivation of $[\cdot,\cdot]_0$;
\item Proposition~\ref{prop:not-additive}: if $\Omega\neq0$, neither
$T^{(1)}\cdots T^{(r)}$ nor $U_{\le1}$ (\texttt{T1} in the development) intertwines $[\cdot,\cdot]_\beta$ with $[\cdot,\cdot]_0$;
for $u\neq v$, the $L_{uv}$-coefficient of the defect $[U_{\le1}L_{uu},U_{\le1}L_{vv}]_0-U_{\le1}[L_{uu},L_{vv}]_\beta$ is
$2\Omega_{uv}$;
\item Remark~\ref{rem:even}: $J-\Omega$ is equivalent to $J$ by
$\Theta=(\mathbb I+\Omega J)^{1/2}$ for every even antisymmetric $\Omega$ whose entries
have $\kappa$-degree~$2$.
\end{enumerate}
\item Controls. For $N=1$ the bracket $[\cdot,\cdot]_0$, $\Gamma_\beta$, $\tilde h$ and $f_\beta$
coincide, under the evident identification of bases and coefficients, with those of the
formalization~\cite{Repo} of~\cite{Aoi26}. In particular, $\tilde h$ is the map $h$ and
$f_\beta$ the primitive of~\cite[Theorem~4.1]{Aoi26}. For $C(2)=\osp(2|2)$ the basis has
$8$ elements, and three values of $\Gamma_\beta$ are evaluated there. They
agree with an exact computation of the bracket in $A_\beta$, made outside the
development.
\item Remark~\ref{rem:d21}. For $D(2,1)=\osp(4|2)$ the basis has $9$ even and $8$ odd
elements, its Killing form vanishes, and it is the only $\osp(N|2n)$ of superdimension
$(9|8)$. The identity $\Gamma_\beta=\delta\tilde h$, Theorem~\ref{thm:triv-single} and
Corollary~\ref{cor:tbeta} are stated there as instances. Every statement of items 1--8
applies there, since the fermionic index is arbitrary.
\end{enumerate}

Two proofs in the development differ from those in the paper. Theorem~\ref{thm:all-primitives} is proved
with a grading element of $\mathfrak{sp}(2n)$ instead of Whitehead's lemma, and \eqref{eq:linear}
follows from the identity $[UX,UY]_0=U([X,Y]_0+\sum_s\kappa_s\Gamma^{(s)}(X,Y))$ on $\mathfrak g_R$ and
Theorem~\ref{thm:multi}, instead of the computation in $A_\beta$.

\subsection{Declarations}\label{subsec:formal-binding}

Names are relative to the namespace
\begin{center}\texttt{OspOddDeformations.OspN},\end{center}
except \texttt{General.intertwine\_map}, which is relative to \texttt{OspOddDeformations}.

{\small
\begin{longtable}{ll}
item & declarations\\\hline
\endhead
1 & \texttt{bracket0Basis\_super\_skew}, \texttt{bracket0Basis\_degree}\\
2 & \texttt{Gamma\_eq\_delta\_hTilde}\\
3 & \texttt{intertwining}, \texttt{TBetaInv\_TBeta}, \texttt{TBeta\_TBetaInv}\\
4 & \texttt{Presented.untwist}, \texttt{Presented.untwist\_b}, \texttt{Presented.untwist\_symm\_b},\\
  & \texttt{Presented.untwist\_a}, \texttt{Presented.untwist\_}$\kappa$, \texttt{Presented.untwist\_parity},\\
  & \texttt{Presented.base\_realization}, \texttt{Presented.untwist\_Xp\_eq},\\
  & \texttt{Presented.untwist\_iotaP}, \texttt{Presented.recovery}, \texttt{Presented.recovery\_delta},\\
  & \texttt{Presented.bracketBeta\_unique}, \texttt{Presented.scommR\_eR}\\
5 & \texttt{Presented.iotaRP\_injective}, \texttt{Presented.iotaRP\_beta\_injective},\\
  & \texttt{Presented.recovery\_unique}\\
6 & \texttt{delta\_ad\_eq\_zero}, \texttt{fBeta\_L}, \texttt{fBeta\_F}, \texttt{fBeta\_S},\\
  & \texttt{Gamma\_eq\_delta\_fBeta}, \texttt{intertwining\_fBeta},\\
  & \texttt{TBetaFInv\_TBetaF}, \texttt{TBetaF\_TBetaFInv}, \texttt{der\_eq\_ad},\\
  & \texttt{primitive\_eq\_hTilde\_sub\_ad}, \texttt{fBeta\_unique\_primitive},\\
  & \texttt{no\_primitive\_vanishing}, \texttt{complex\_parameters}, \texttt{TBetaC\_bijective},\\
  & \texttt{General.intertwine\_map}, \texttt{hTilde\_Fof\_ne\_zero}, \texttt{fBeta\_Sof\_ne\_zero}\\
7 & \texttt{rho\_ospCond}, \texttt{rho\_bracket},
    \texttt{rhoBasis\_linearIndependent}, \texttt{rho\_span},\\
  & \texttt{gram\_ordOdd}, \texttt{gram\_ordEven}, \texttt{ordOdd\_bijective},
    \texttt{ordEven\_bijective}\\
8(a) & \texttt{Presented.untwistR}, \texttt{Presented.untwistR\_b}, \texttt{Presented.untwistR\_symm\_b},\\
  & \texttt{Presented.untwistR\_a}, \texttt{Presented.untwistR\_}$\kappa$, \texttt{Presented.untwistR\_parity},\\
  & \texttt{Multi.thetaSeries\_spec}, \texttt{Multi.thetaSeries\_eq\_thetaHalf}\\
8(b) & \texttt{Presented.untwistR\_XpM}, \texttt{Presented.untwistR\_iotaM},\\
  & \texttt{Multi.UM\_bijective}, \texttt{Multi.UM\_thetaSeries\_bijective}\\
8(c) & \texttt{Presented.recoveryM}, \texttt{Presented.iotaM\_beta\_injective},\\
  & \texttt{Presented.bracketBetaM\_unique}, \texttt{Presented.scommM\_homog},\\
  & \texttt{Multi.intertwining\_multi},\\
  & \texttt{Multi.intertwining\_thetaSeries}, \texttt{Multi.bracketBetaM\_eq\_conj}\\
8(d) & \texttt{Multi.VChi\_comp}, \texttt{Multi.TS\_comm}, \texttt{Multi.TProd\_eq},
    \texttt{Multi.ThetaStar\_comp},\\
  & \texttt{Multi.UM\_eq\_ThetaStar\_TProd}, \texttt{Multi.vChi\_bracket},
    \texttt{Multi.thetaStar\_bracket},\\
  & \texttt{Multi.ThetaStar\_not\_aut}\\
8(e) & \texttt{GroupB.jacobi\_eP}, \texttt{GroupB.bracket0M\_derivation},
    \texttt{GroupB.exp\_adM\_aut},\\
  & \texttt{GroupB.ThetaStar\_ne\_exp\_ad}, \texttt{Multi.ThetaStar\_sub\_one\_not\_derivation}\\
8(f) & \texttt{Multi.TProd\_not\_intertwine}, \texttt{Multi.T1\_not\_intertwine},
    \texttt{Multi.T1\_defect\_Luv}\\
8(g) & \texttt{Multi.thetaG\_spec}\\
9 & \texttt{control\_N1\_bracket}, \texttt{control\_N1\_gamma},\\
  & \texttt{control\_N1\_hTilde}, \texttt{control\_N1\_fBeta},\\
  & \texttt{evenIdx\_basis\_card\_1\_1}, \texttt{evenC2\_GammaBasis\_SF},\\
  & \texttt{evenC2\_GammaBasis\_FF}, \texttt{evenC2\_GammaBasis\_LS}\\
10 & \texttt{d21\_even\_card}, \texttt{d21\_odd\_card}, \texttt{GroupB.d21\_killing\_zero},
     \texttt{GroupB.sdim\_9\_8\_iff},\\
  & \texttt{d21\_Gamma\_eq\_delta\_hTilde}, \texttt{d21\_intertwining}, \texttt{d21\_intertwining\_fBeta}
\end{longtable}}

\subsection{What is not formalized}\label{subsec:not-formalized}

The following items delimit what has been checked.

\begin{enumerate}
\item \textbf{The Fock module.} Remark~\ref{rem:fock} is not formalized.
\item \textbf{The comparison with~\cite{FSS00}.} The oscillator realizations of~\cite{FSS00}
and their comparison in Section~\ref{sec:special-remarks}, and Remark~\ref{rem:d21} beyond
item~10 (the Killing form of every $D(n+1,n)$, and which $D(2,1;\alpha)$ are isomorphic to
$D(2,1)$), are not formalized.
\item \textbf{Two points of Section~\ref{sec:multi}.} That for $r\leq3$ the map
$\Theta_*-\mathrm{id}$ has $\kappa$-degree $2$ and its square vanishes (the text after
Proposition~\ref{prop:factor}; item~8(e) does not use it), and
Remark~\ref{rem:even-generated}, are not formalized.
\item \textbf{Other remarks.} Remark~\ref{rem:central}, and Remark~\ref{rem:even} beyond
item~8(g), are not formalized.
\end{enumerate}

\subsection{The build and its audit}\label{subsec:formal-repro}

The development uses Lean \texttt{v4.35.0} and mathlib \texttt{v4.35.0}, at the commit
\begin{center}\texttt{56a005855ea17e0b0958f73ec801ae140b0a1679},\end{center}
both pinned in the accompanying sources. Its audit file prints the axiom dependencies of
the $324$ declarations it lists, which include all those above. Each depends at most on
\texttt{propext}, \texttt{Classical.choice} and \texttt{Quot.sound}.
No declaration depends on \texttt{sorryAx}, and no additional axiom is
introduced.

\section*{Use of generative AI}
Generative-AI tools assisted with mathematical discussion, checking, drafting, and the Lean development.
The argument and its classical theorem inputs are stated explicitly above.
Appendix~\ref{sec:formalization} records a machine-checked formalization of the
statements listed there, in the Lean~4 proof assistant; what is not formalized is stated
in that appendix. Exact computations were used as checks; no proof depends on them.
The author reviewed and edited all generated material, and takes full responsibility for the
contents of this paper, irrespective of how they were generated.

\end{document}